# ON LOCAL U-STATISTIC PROCESSES AND THE ESTIMATION OF DENSITIES OF FUNCTIONS OF SEVERAL SAMPLE VARIABLES

BY EVARIST GINÉ[1] AND DAVID M. MASON[2]

*University of Connecticut and University of Delaware*

A notion of local $U$-statistic process is introduced and central limit theorems in various norms are obtained for it. This involves the development of several inequalities for $U$-processes that may be useful in other contexts. This local $U$-statistic process is based on an estimator of the density of a function of several sample variables proposed by Frees [*J. Amer. Statist. Assoc.* **89** (1994) 517–525] and, as a consequence, uniform in bandwidth central limit theorems in the sup and in the $L_p$ norms are obtained for these estimators.

**1. Introduction.** Let $X, X_1, X_2, \dots$ be i.i.d. random variables taking values in $\mathbf{R}$, with common density function $f$ and consider the kernel density estimator of $f$ defined for $t \in \mathbf{R}$,

$$f_n(t, h_n) = (nh_n)^{-1} \sum_{i=1}^{n} K(h_n^{-1}(t - X_i)) =: n^{-1} \sum_{i=1}^{n} K_{h_n}(t - X_i), \tag{1.1}$$

where $\{h_n\}_{n \geq 1}$ is a sequence of positive constants converging to zero at the rate $nh_n \to \infty$ and the kernel $K$ is an integrable (real) function of bounded variation satisfying $\int_{\mathbf{R}} K(x)\,dx = 1$ (Parzen [27]). It is easy to prove that, subject to smoothness conditions on $f$, for each $t \in \mathbf{R}$,

$$\sqrt{h_n} u_n(t) := \sqrt{nh_n}\{f_n(t, h_n) - Ef_n(t, h_n)\} \to_d N(0, \|K\|_2^2 f(t)),$$

whereas for any choice of $t_1 \neq t_2$ the random variables $\sqrt{h_n} u_n(t_1)$ and $\sqrt{h_n} u_n(t_2)$ are asymptotically independent. This means that $\sqrt{h_n} u_n$ cannot converge weakly to a continuous bounded process on any nontrivial subinterval of $\mathbf{R}$.

[1]Supported in part by NSA Grant H98230-04-1-0075.
[2]Supported in part by NSA Grant MDA904-02-1-0034 and NSF Grant DMS-02-03865.

Since the bias $Ef_n(t,h_n) - f(t)$ can always be dealt with under the usual conditions on $K$ and $f$, this tells us that $f_n(t,h_n)$ estimates the density $f$ at a much slower rate than $n^{-1/2}$. On the other hand, Frees [15] discovered the perhaps surprising fact that the densities $f_g(t)$ of some symmetric real functions $g(X_1,\dots,X_m)$ of $m>1$ i.i.d. random variables can be estimated at each fixed $t$ at the rate $n^{-1/2}$, using the $U$-statistic estimator

$$\frac{(n-m)!}{n!}\sum_{\mathbf{i}\in I_n^m} K_{h_n}(t-g(X_{i_1},\dots,X_{i_m})),$$

where $\mathbf{i}=(i_1,\dots,i_m)$ and $I_n^m=\{(i_1,\dots,i_m): 1\le i_j\le n, i_j\ne i_k \text{ if } j\ne k\}$. Schick and Wefelmeyer [32] consider the special case

$$g(X_1,\dots,X_m)=\sum_{i=1}^m u_i(X_i),$$

which however is not necessarily a symmetric function. Using convolution kernels, they also obtain the (in probability) rate of $n^{-1/2}$ for the sup norm and the $L_1$ norm measure of the discrepancy between the estimator and the density (actually, they obtain limit theorems in distribution). These results require smoothness conditions on the kernel, the density and certain conditional densities.

Our aim is to extend the Schick and Wefelmeyer results to the Frees framework and with $g$ not necessarily symmetric. Moreover, we frame our results uniform in bandwidth in the sup norm and the $L_p$ norm with $p\ge 1$ so that they can be used with adaptive bandwidth estimators. Also, we do this not necessarily in one dimension, but in $\mathbf{R}^d$. These extensions substantially increase the scope of applicability of the results and, as a consequence, we can show how our results recover, extend and/or improve upon, previous work by several authors. We shall not discuss removing the bias in general but only in some examples and very briefly. The bias is not probabilistic and can always be treated by adding enough smoothness to the kernel and the density (see, e.g., the two references just cited and Ahmad and Fan [1]).

We shall begin by generalizing the setup in Giné, Mason and Zaitsev [21] (concretely, that of Examples 1.2 and 1.3 there) and Mason [25] to $U$-statistics. Throughout this paper, we let $X, X_i$, $i\in\mathbf{N}$, be i.i.d. random variables taking values in a measurable space $(S,\mathcal{S})$; let $g: S^m\mapsto\mathbf{R}^d$, $1\le d<\infty$, be a measurable function; let $K:\mathbf{R}^d\mapsto\mathbf{R}$ be an integrable measurable function that integrates to 1 (a "kernel"); and let $0<a\le b<\infty$. Then, for $t\in\mathbf{R}^d$ and $\lambda\in[a,b]$, we introduce the *local $U$-statistic*

$$U_n(t,\lambda) := \frac{(n-m)!}{n!}\sum_{\mathbf{i}\in I_n^m} K_{\lambda h_n}(t-g(X_{i_1},\dots,X_{i_m})), \tag{1.2}$$

where, here and elsewhere in this paper, for functions $H:\mathbf{R}^d \mapsto \mathbf{R}$ and $h>0$, the notation

$$H_h(t) := h^{-1}H(t/h^{1/d}), \qquad t\in\mathbf{R}^d, \tag{1.3}$$

is in force. The term "local" just reflects the fact that the $U$-process (1.2) is of a special kind, namely a convolution of an approximate identity $K_{\lambda h_n}$ with an "empirical measure," in this case $((n-m)!/n!)\sum_{I_n^m}\delta_{g(X_{i_1},\dots,X_{i_m})}$, and therefore, for each value of $t$, the largest contributions to the statistic come from the values $g(X_{i_1},\dots,X_{i_m})$ closest to $t$.

Special cases of $U_n(t,\lambda)$, when $X_1,\dots,X_n$ are i.i.d. $\mathbf{R}^d$ valued, include the interpoint distance studied by Jammalamadaka and Janson [23], with $g(x,y)=|x-y|$,

$$h_nU_n(0,1)=\frac{1}{n(n-1)}\sum_{\mathbf{i}\in I_n^2} I\{|X_{i_1}-X_{i_2}|\le h_n\}$$

and the related short distance process studied by Eastwood and Horváth [12],

$$\frac{1}{n(n-1)}\sum_{\mathbf{i}\in I_n^2} I\{\mathbf{d}(X_{i_1},X_{i_2})\le \lambda h_n\}, \qquad 0\le\lambda\le 1,$$

where $\mathbf{d}$ is a distance on $\mathbf{R}^d$, as well as the $U$-statistic estimator of the density of the sum $X_1+\cdots+X_m$,

$$\frac{(n-m)!}{n!}\sum_{\mathbf{i}\in I_n^m} K_{h_n}(t-X_{i_1}-\cdots-X_{i_m}).$$

Our goal is to obtain central limit theorems for the following *local $U$-statistic process* formed from $U_n(t,\lambda)$:

$$u_{n,\lambda}(t) := \sqrt{n}\{U_n(t,\lambda)-EK_{\lambda h_n}(t-g(X_1,\dots,X_m))\}, \qquad t\in\mathbf{R}^d, \lambda\in[a,b]. \tag{1.4}$$

The case when $m=1$ is a special case of the local empirical process studied in Mason [25]. We shall confine our attention to the case $m\ge 2$, which we shall soon see has a radically different asymptotic behavior than the case $m=1$. Occasionally, we may restrict the process to $t\in D\subset\mathbf{R}^d$, where $D$ may even consist of a single point.

The limit theorems to be obtained for these processes will be in the sup and in the $L_p$ norms, $1\le p<\infty$, uniformly in $\lambda\in[a,b]$ (precise definitions in the next section). The reason these results will be true will be essentially the same as in Frees [15], namely: (1) the process (1.4) is equivalent to its linear part, that is, all the terms in its Hoeffding decomposition of order higher than one tend to zero in the appropriate way; (2) that this linearization, which is a smoothed empirical process, is equivalent to the empirical process without smoothing, hence *independent* of $\lambda$ and $h_n$; (3) that this empirical process, under appropriate and quite weak conditions, satisfies the central limit theorem in the sup norm or in the $L_p$ norms.

The devil is in the details, and extending the Frees result from one point $t \in \mathbf{R}$, to uniformity in $t \in \mathbf{R}^d$ and in $\lambda \in [a,b]$, makes for very different proofs and requires a substantial amount of technique, some of it new.

The statements of the central limit theorems, along with examples, are collected in Section 2 and all the proofs are postponed to Section 3. In the process of establishing our results we shall develop some tools that should be of separate interest. Among them include tight bounds for the absolute moment of the supremum of the $U$-statistic process under a uniform covering number bound generalizing a similar bound obtained by Einmahl and Mason [13] and Giné and Koltchinskii [17] (see also Giné and Guillou [16]) for the usual empirical process. All of these results were motivated by Proposition 6.2 of Talagrand [35], which is an expectation bound for VC classes of sets, where VC stands for Vapnik and Červonenkis. We also obtain moment and exponential inequalities for Banach space valued $U$-statistics which, although not necessarily optimal, are very easy to apply and are well adapted to the problems treated in this article. In a sequel to this paper (Giné and Mason [20]), we derive the corresponding functional laws of the logarithm for $u_{n,\lambda}$.

**2. Main results and examples.** As is usual when dealing with empirical processes, we take $((S \times T)^{\mathbf{N}}, (\mathcal{S} \otimes \mathcal{T})^{\mathbf{N}}, \Pr)$ as the underlying probability space, where $(S, \mathcal{S})$ is a measurable space, $T = \{-1, 1\}$, $\mathcal{T}$ is the family of all subsets of $T$, and $\Pr = P^{\mathbf{N}} \times (P')^{\mathbf{N}}$, $P$ a probability measure on $(S, \mathcal{S})$ and $P'$ the uniform distribution on $T$. Then the random variables $X_i$ are the projections $(S \times T)^{\mathbf{N}} \mapsto S$, $X_i(s_1, t_1, s_2, t_2, \ldots) = s_i$, which are i.i.d. with law $P$. We will occasionally use the random variables $\varepsilon_i(s_1, t_1, s_2, t_2, \ldots) = t_i$, which obviously satisfy $\Pr\{\varepsilon_i = 1\} = \Pr\{\varepsilon_i = -1\} = 1/2$. Note that the random variables $\{X_i, \varepsilon_j : i, j \in \mathbf{N}\}$ are independent. The variables $\varepsilon_i$ are often called Rademacher variables. Sometimes we will write $X$ for $X_1$.

All the asymptotic results on the process $u_{n,\lambda}$ in this article require the following key assumption.

(CD) For each $i = 1, \ldots, m$, the random variable $g(X_1, \ldots, X_m)$, conditionally on $X_i = x$, $x \in S$, has a density $\overline{f}_i(t, x)$, $t \in \mathbf{R}^d$, which is jointly measurable in $t$ and $x$.

Note that, setting

$$\overline{f} = \sum_{i=1}^{m} \overline{f}_i, \tag{2.1}$$

the function

(2.2) $$f_g(t) := E[\overline{f}(t,X)]/m, \qquad t \in \mathbf{R}^d,$$

defines a density for the random variable $g(X_1,\dots,X_m)$. Another condition that we will require for different values of $p \in [1,\infty]$ is

(CD$p$) $$\|\overline{f}(\cdot,x)\|_p < \infty \qquad \text{for all } x \in S \text{ and } \|f_g\|_p < \infty.$$

Here and elsewhere, $\|\cdot\|_p$ denotes the $L_p(\mathbf{R}^d)$ norm for $1 \le p < \infty$ and the sup norm on $\mathbf{R}^d$ for $p = \infty$.

The asymptotics of the processes $u_{n,\lambda}$ will turn out to be equivalent to that of the processes

(2.3) $$\overline{\nu}_n(t) = \frac{1}{\sqrt{n}} \sum_{i=1}^n (\overline{f}(t,X_i) - E[\overline{f}(t,X_i)]), \qquad t \in \mathbf{R}^d,$$

which are empirical processes. With some abuse of notation, we say that

(2.4) $$\overline{\nu}_n \text{ converges in law in } L_\infty(\mathbf{R}^d)$$

if condition (CD$\infty$) holds and the class of functions $\{\overline{f}(t,\cdot) : t \in \mathbf{R}^d\}$ is $P$-Donsker. We refer to Dudley [10] or van der Vaart and Wellner [37] for the definition of Donsker classes of functions. The condition (CD$\infty$) is equivalent to the maps $t \mapsto \overline{f}(t,x)$ and $t \mapsto E\overline{f}(t,X)$ being in $\ell^\infty(\mathbf{R}^d)$, the abuse of notation consists in replacing this last space by $L_\infty(\mathbf{R}^d)$, which usually means something else, and the precise meaning of (2.4) is that $\lim_{n\to\infty} EH(\nu_n)^* = EH(G)$ for every $H : \ell_\infty(\mathbf{R}^d) \mapsto \mathbf{R}$ bounded and continuous, where $G = \{G(t) : t \in \mathbf{R}^d\}$ is the centered Gaussian process with the covariance of $\overline{f}(\cdot,X)$, more precisely a sample continuous version, and $H(\overline{\nu}_n)^*$ is the (a.s.) smallest measurable function larger than or equal to $H(\overline{\nu}_n)$. Likewise, for $1 \le p < \infty$, we say that

(2.5) $$\overline{\nu}_n \text{ converges in law in } L_p(\mathbf{R}^d)$$

if the condition (CD$p$) holds and the $L_p(\mathbf{R}^d)$-valued random variable $\overline{f}(\cdot,X)$ satisfies the central limit theorem in this space, that is, there is a centered Gaussian process $G$ with the same covariance as $\overline{f}(\cdot,X)$ and with sample paths in $L_p(\mathbf{R}^d)$, such that for every $H : L_p(\mathbf{R}^d) \mapsto \mathbf{R}$ bounded and continuous,

$$\lim_{n\to\infty} EH(\overline{\nu}_n) = EH(G).$$

Note that in each case, $\overline{\nu}_n(\cdot)$ is a random variable that takes values in $L_p(\mathbf{R}^d)$, $1 \le p \le \infty$, although in the case $p = \infty$, with abuse of notation: $\overline{\nu}_n(\cdot)$ is really in $\ell^\infty(\mathbf{R}^d)$, the space of bounded functions on $\mathbf{R}^d$, and $\overline{\nu}_n$ is not necessarily measurable.

The following definition describes the type of central limit theorem we will prove for the process $u_{n,\lambda}$.

DEFINITION 1. Let $1 \le p \le \infty$ and assume (CD) and (CD$p$). The processes $u_{n,\lambda}$ converge weakly in $L_p(\mathbf{R}^d)$, uniformly in $a \le \lambda \le b$, to the centered Gaussian process $G$ with the same covariance as $\overline{f}(\cdot, X)$ if

$$\sup_{\lambda \in [a,b]} \|u_{n,\lambda} - \overline{\nu}_n\|_p \to 0 \qquad \text{in pr}^*$$

and $\overline{\nu}_n$ converges weakly in $L_p(\mathbf{R}^d)$ in the sense of (2.4) for $p = \infty$ and (2.5) for $p < \infty$.

Convergence in pr$^*$ means convergence in probability of the measurable envelopes (the smallest dominating measurable functions).

Note that in the case $p = \infty$, convergence of $u_{n,\lambda}$ to $G$ in the sense of this definition implies that the processes $(t, \lambda) \mapsto u_{n,\lambda}(t)$ converge in law in $\ell^\infty(\mathbf{R}^d \times [a,b])$ to the Gaussian process $G$ (de la Peña and Giné [7], Dudley [10], or van der Vaart and Wellner [37] for this type of convergence). However, our notion gives more. In fact if, for $l = 1, \dots, N$, we have functions $g_l : S^{m_l} \mapsto \mathbf{R}^{d_l}$ such that the processes $u_{n,\lambda}^{(l)}$ corresponding to $g = g_l$ converge weakly in $L_\infty(\mathbf{R}^{d_l})$, uniformly in $a_l \le \lambda \le b_l$, to a centered Gaussian process $G_l$ with the same covariance as $F_l(\cdot, X)$, where $F_l$ is the $\overline{f}$ corresponding to $g = g_l$, then it is easy to conclude, using the obvious multivariate extension of Definition 1, that the vector-valued processes $\overrightarrow{u}_{n,\lambda_1,\dots,\lambda_N}$ defined by

$$\overrightarrow{u}_{n,\lambda_1,\dots,\lambda_N}(t_1, \dots, t_N) = (u_{n,\lambda_1}^{(1)}(t_1), \dots, u_{n,\lambda_N}^{(N)}(t_N)), \qquad t_1 \in \mathbf{R}^{d_1}, \dots, t_N \in \mathbf{R}^{d_N}, \tag{2.6}$$

converge weakly in $L_\infty(\mathbf{R}^{d_1}) \times \cdots \times L_\infty(\mathbf{R}^{d_N})$ uniformly in $a_l \le \lambda_l \le b_l$, $l = 1, \dots, N$, to the centered vector-valued Gaussian process defined on $\mathbf{R}^{d_1} \times \cdots \times \mathbf{R}^{d_N}$,

$$G(t_1, \dots, t_N) = (G_1(t_1), \dots, G_N(t_N)), \tag{2.7}$$

with the same covariance/cross covariance matrix as $F(\cdot, X)$, where

$$F(t_1, \dots, t_N, X) = (F_1(t_1, X), \dots, F_N(t_N, X)), \qquad t_1 \in \mathbf{R}^{d_1}, \dots, t_N \in \mathbf{R}^{d_N}. \tag{2.8}$$

2.1. *Central limit theorems.* We still require another definition. We say that a class of measurable functions $\mathcal{F}$ defined on a measurable space $(S, \mathcal{S})$ is VC-type (VC for Vapnik and Červonenkis) with respect to an envelope $F$ (meaning a measurable function $F$ such that $|f| \le F$ for all $f \in \mathcal{F}$) if the covering number $N(\mathcal{F}, L_2(Q), \varepsilon)$, defined as the smallest number of $L_2(Q)$ open balls of radius $\varepsilon$ required to cover $\mathcal{F}$, satisfies

$$N(\mathcal{F}, L_2(Q), \varepsilon) \le \left( \frac{A\|F\|_{L_2(Q)}}{\varepsilon} \right)^v, \qquad 0 < \varepsilon \le 2\|F\|_{L_2(Q)}, \tag{2.9}$$

for some $A \geq 3$ and $v \geq 1$, for every probability measure $Q$ on $\mathcal{S}$ for which $Q(F^2) < \infty$. If (2.9) holds for $\mathcal{F}$, then we say that the VC class $\mathcal{F}$ admits the characteristics $A$ and $v$.

THEOREM 1. *Let $p = \infty$, let* (CD) *and* (CD$\infty$) *hold, and let $K$ be bounded. Assume:*

(a) *each of the classes*

$$\mathcal{K}_n := \{K(h^{-1/d}(y - \cdot)) : y \in \mathbf{R}^d, ah_n \leq h \leq bh_n\} \tag{2.10}$$

*is VC-type for a bounded envelope $F_n$ and all the classes $\mathcal{K}_n$ admit the same characteristics $A$ and $v$;*

(b) *the density $f_g$ of $g(X_1, \ldots, X_m)$ is bounded and the class of functions $\mathcal{F} := \{\overline{f}(t, \cdot) : t \in \mathbf{R}^d\}$ is $P$-Donsker and the identity map*

$$(\mathbf{R}^d, |\cdot|) \mapsto (\mathbf{R}^d, \rho)$$

*is uniformly continuous, where $\rho^2(u, v) = \mathrm{Var}(\overline{f}(u, X) - \overline{f}(v, X))$;*

(c) *$h_n \to 0$ and $nh_n/(1 \vee \log(A\|F_n\|_{L_2(f_g)}/\sqrt{h_n}))^2 \to \infty$. The processes $u_{n,\lambda}$ then converge weakly in $L_\infty(\mathbf{R}^d)$, uniformly in $a \leq \lambda \leq b$, to the centered Gaussian process with the same covariance as $\overline{f}(\cdot, X)$.*

REMARK 1. Theorem 1 has obvious applications to the construction of confidence bands for $f_g$. It is formulated uniformly in $a \leq \lambda \leq b$ so as to allow the possibility for $\lambda$ to be replaced by an estimator $\widehat{\lambda}_n$. Suppose that $\widehat{\lambda}_n = \lambda_n(X_1, \ldots, X_n)$ is an adaptive bandwidth selector such that for all $0 < \varepsilon < 1$ there exist $0 < c < d < \infty$ for which for all large enough $n$

$$\Pr\{c < \widehat{\lambda}_n < d\} \geq 1 - \varepsilon. \tag{2.11}$$

Then if assumption (a) of Theorem 1 holds for any choice of $0 < a < b < \infty$ we can immediately conclude from (2.11) that the processes $u_{n,\widehat{\lambda}_n}$ converge weakly in $L_\infty(\mathbf{R}^d)$ to the centered Gaussian process with the same covariance as $\overline{f}(\cdot, X)$. The analogous remark holds for Theorems 2 and 3 below. For a thorough discussion of bandwidth estimators that satisfy (2.11) refer to Deheuvels and Mason [5].

REMARK 2. We note that if the classes of functions $\mathcal{F}_i = \{\overline{f}_i(t, \cdot) : t \in \mathbf{R}^d\}$, $i = 1, \ldots, m$, are $P$-Donsker, then so is $\mathcal{F}$ given in (b) (e.g., Theorem 2.10.6 in van der Vaart and Wellner [37] applied to the function $\phi(t_1, \ldots, t_m) = \sum_{i=1}^m t_i/\sqrt{m}$). A sufficient condition for the class of functions $\mathcal{F}$ being $P$-Donsker is that it be of VC-type, since the function $\overline{f}$ being jointly measurable already ensures that this class of functions is measurable (e.g., Definition 2.3.3, Example 2.3.5 and Theorem 2.5.2 (Pollard's CLT) in van der Vaart and Wellner [37], or Dudley [10], Theorem 6.3.1). These two references contain many examples of classes of functions which are of VC-type. We single out one of them which is particularly convenient for condition (a) above. Let $\Psi : \mathbf{R} \mapsto \mathbf{R}$ be a function of bounded variation on $\mathbf{R}$ ($\Psi$ is the difference of two bounded nondecreasing functions). The proof of Lemma 22 in Nolan and Pollard [26] shows that if

$$K(x) = \Psi(p(x)), \qquad x \in \mathbf{R}^d, \tag{2.12}$$

where $p$ is either a real polynomial on $\mathbf{R}^d$ or the $\alpha$th power of the absolute value of a real polynomial on $\mathbf{R}^d$, $\alpha > 0$, then the class of functions

$$\mathcal{K} = \{K(\gamma^{-1}(t - \cdot)) : t \in \mathbf{R}^d, \gamma > 0\}$$

is VC-type. Moreover, since the function $K(\gamma^{-1}(t-x))$ is jointly measurable, this class is also measurable. Most kernels of interest satisfy condition (2.12), and therefore, condition (a) in Theorem 1.

Sometimes one is interested only in weak convergence in $L_\infty(D)$ uniformly in $a \le \lambda \le b$, where $D$ is a subset of $\mathbf{R}^d$ that may even consist of a single point. For instance, this is the case for the interpoint distance process, where $K(x) = I\{|x| \le 1\}$, $S = \mathbf{R}^d$ and $g(x, y) = x - y$. The following, which is related to a result of Eastwood and Horváth **(year?)**, will be a corollary to the proof of the above theorem. We say that a subset $D$ of $\mathbf{R}^d$ is star-shaped about 0 if $x \in D$ implies $\lambda x \in D$ for all $0 \le \lambda \le 1$.

COROLLARY 1. *Let $D$ be a measurable bounded subset of $\mathbf{R}^d$ star-shaped about 0 and with $\mathrm{Vol}(D) \neq 0$, and set $K = I_D/\mathrm{Vol}(D)$. Assume that $g(X_1, \dots, X_m)$ has a bounded density $f_g$, that $E[\overline{f}(0, X)]^2 < \infty$ and that for all $\varepsilon > 0$,*

$$\lim_{\delta \to 0} \limsup_{n \to \infty} \Pr{}^* \left\{ \sup_{|u| \le \delta} |\overline{\nu}_n(u) - \overline{\nu}_n(0)| > \varepsilon \right\} = 0. \tag{2.13}$$

*Let $h_n \to 0$ be such that $nh_n \to \infty$. Then, the processes $\lambda u_{n,\lambda}(0)$ converge weakly in $\ell^\infty((0,1])$ to the process $\lambda \sigma Z$, $0 \le \lambda \le 1$, where $Z$ is standard normal and $\sigma^2 = \mathrm{Var}(\overline{f}(0, X))$.*

It makes sense to define $\lambda u_{n,\lambda}(0)$ as zero for $\lambda = 0$. With this convention, in force from here on, we have weak convergence in $\ell^\infty([0,1])$ in this corollary.

Condition (2.13) is satisfied, for instance, if the class $\mathcal{F}_\delta := \{\overline{f}(t, \cdot) : |t| \le \delta\}$ is $P$-Donsker for some $\delta > 0$, a condition weaker than the class $\mathcal{F}$ being $P$-Donsker.

Next we state the central limit theorems in the $L_p$ norm, uniform on $a \le \lambda \le b$. We need to recall the definition of Young moduli of exponential type. As in de la Peña and Giné [7], page 188, $\Psi_1(x) := e^x - 1$, but if $\alpha < 1$, since $e^{x^\alpha}$ is only convex for $x \geq x_\alpha := ((1-\alpha)/\alpha)^{1/\alpha}$, we take as a function $\Psi_\alpha$ that is 0 at 0, convex and increasing, and of the order of $e^{x^\alpha}$ for $x$ large the following:

$$\Psi_\alpha(x) := \tau_\alpha(x) - \alpha \exp((1-\alpha)/\alpha), \tag{2.14}$$

where $\tau_\alpha(x)$ equals $\exp(x^\alpha)$ if $x \geq x_\alpha$, and equals the tangent line to the function $y = \exp(x^\alpha)$ at $x = x_\alpha$ for $0 \leq x \leq x_\alpha$. We also recall that then, for any nonnegative random variable $\xi$ for which $Ee^{(\xi/a)^\alpha} < \infty$ for some $a > 0$,

$$\|\xi\|_{\Psi_\alpha} := \inf\{c : E\Psi_\alpha(\xi/c) \leq 1\}. \tag{2.15}$$

This is a (pseudo)norm and it dominates, up to constants that depend only on $\alpha$ and $p$, all the $L_p$ (pseudo)norms. Simple standard computations show that

$$\Pr\{\xi > x\} \leq b \exp\{-(x/a)^\alpha\} \qquad \text{for all } x > 0 \quad \Longrightarrow \quad \|\xi\|_{\Psi_\alpha} \leq Ca \tag{2.16}$$

for a constant $C$ that depends only on $\alpha$ and $b$. Note also that $\Psi_\alpha^{-1}(u)$ is a constant times $u$ for $0 \leq u \leq \Psi_\alpha(x_\alpha)$ and it is the $1/\alpha$-th power of the logarithm of $u + \alpha \exp((1-\alpha)/\alpha)$ for $u \geq \Psi_\alpha(x_\alpha)$.

Given a kernel $K$ and $0 < a \leq b < \infty$, define

$$\mathcal{K}_{[a,b]} := \{K_h : h \in [a,b]\} \cup \{0\}. \tag{2.17}$$

For any two functions $f$ and $g$ in $L_p$, set

$$d_p^p(f,g) = \int_{\mathbf{R}^d} |f(t) - g(t)|^p \, dt. \tag{2.18}$$

THEOREM 2. *Let $2 \leq p < \infty$ and let (CD) and (CD$p$) hold. Assume:*

(a) *the kernel $K$ is in $L_p(\mathbf{R}^d)$ and*

$$\int_0^1 \Psi_{2/m}^{-1}(N(\mathcal{K}_{[a,b]}, d_p, \varepsilon)) \, d\varepsilon < \infty; \tag{2.19}$$

(b) *the sum of conditional densities $\overline{f}$ defined in (2.1) satisfies*

$$\begin{gathered} \lim_{t\to\infty} t^2 \Pr\{\|\overline{f}(\cdot, X) - E\overline{f}(\cdot, X)\|_p > t\} = 0 \\ \text{and} \quad \int_{\mathbf{R}^d} [E(\overline{f}(u,X) - E\overline{f}(u,X))^2]^{p/2} \, du < \infty; \end{gathered} \tag{2.20}$$

(c) *$h_n \to 0$ and $nh_n^{2(p-1)/p} \to \infty$.*

*Then, the processes $u_{n,\lambda}$ converge weakly in $L_p(\mathbf{R}^d)$, uniformly in $a \leq \lambda \leq b$, to the centered Gaussian process with the same covariance as $\overline{f}(\cdot, X)$.*

The two conditions in (2.20) are implied by the stronger but sometimes more convenient conditions $E\|\overline{f}_i(\cdot,X)-E\overline{f}_i(\cdot,X)\|_p^2<\infty$ for $i=1,\dots,m$. This follows immediately from Minkowski's inequality for integrals (e.g., Folland [14], page 194). Also, note that for $p=2$ the first condition in (2.20) is superfluous.

We should remark here that the conditions (2.20) are precisely the necessary and sufficient conditions for the $L_p(\mathbf{R}^d)$-valued random variable $\overline{f}(\cdot,X)$ to satisfy the central limit theorem, that is, for the processes $\overline{\nu}_n$ to converge in law to $G$ in $L_p(\mathbf{R}^d)$, $p\ge 2$ (Pisier and Zinn [29]; see also Araujo and Giné [2], pages 206–207).

For $1\le p<2$ we need $K^2$ and $f_g$ to satisfy certain moment assumptions, and for this it will be convenient to have the following notation: for $s>0$, $p\ge 1$, we define the Borel measure $\mu_s$ on $\mathbf{R}^d$, $L_p(\mu_s)$ and $\tilde{d}_{p,s}$ as

$$\begin{aligned} d\mu_s(t)&=(1+|t|)^s\,dt, \qquad L_p(\mu_s)=L_p(\mathbf{R}^d,\mathcal{B},\mu_s),\\ \tilde{d}_{p,s}(f,g)&=\|f-g\|_{L_p(\mu_s)}, \end{aligned}\tag{2.21}$$

the latter for functions $f,g\in L_p(\mu_s)$.

THEOREM 3. *Let $1\le p<2$ and let (CD) and (CD$p$) hold. Assume $K^2$ and $f_g$ are in $L_1(\mu_s)$ for some $s>d(2-p)/p$. Assume also:*

(a) *$K$ is in $L_p(\mathbf{R}^d)$ and*

$$\int_0^1 \Psi_{1/m}^{-1}(N(\mathcal{K}_{[a,b]},d_p\vee\tilde{d}_{2,s},\varepsilon))\,d\varepsilon<\infty;\tag{2.22}$$

(b) *the sum $\overline{f}$ of conditional densities satisfies*

$$\int_{\mathbf{R}^d}[E(\overline{f}(t,X)-E\overline{f}(t,X))^2]^{p/2}\,dt<\infty;\tag{2.23}$$

(c) *$h_n\to 0$ and $nh_n\to\infty$.*

*Then, the processes $u_{n,\lambda}$ converge weakly in $L_p(\mathbf{R}^d)$, uniformly in $a\le\lambda\le b$, to the centered Gaussian process with the same covariance as $\overline{f}(\cdot,X)$.*

For $p=1$, condition (2.23) is equivalent to $\int[E\overline{f}^2(t,X)]^{1/2}\,dt<\infty$.

As above, we should also mention that condition (2.23) is precisely the necessary and sufficient condition for the $L_p(\mathbf{R}^d)$-valued random variable $\overline{f}(\cdot,X)$ to satisfy the central limit theorem, that is, for the processes $\overline{\nu}_n$ to converge in law to $G$ in $L_p(\mathbf{R}^d)$, $1\le p\le 2$ (Vakhania [36] and Jain [22]; see also Araujo and Giné [2], pages 206–207).

REMARK 3. All garden variety kernels satisfy the entropy assumptions of Theorems 2 and 3. In fact, many such kernels are of the form $K(x)=\Psi(|x|)$, where $\Psi$ is a function of bounded variation defined on $[0,\infty)$. To see this, assume that for bounded nondecreasing functions $M$ and $N$ we can write $\Psi = P - N$. Further, assume that

$$\int_0^\infty r^{d-1}|M(r)|\,dr = M_0 < \infty \quad \text{and} \quad \int_0^\infty r^{d-1}|N(r)|\,dr = N_0 < \infty.$$

Define the class of functions $\mathcal{P} = \{r \mapsto r^{d-1}\Psi(r\lambda^{1/d}) : \lambda \geq 1\}$. We claim that this class satisfies for some $A_1 > 0$,

$$N(\mathcal{P}, d_1, \varepsilon) \leq A_1 \varepsilon^{-1}, \qquad 0 < \varepsilon \leq 1. \tag{2.24}$$

To verify this choose $1 \leq \lambda \leq \mu$. We see that

$$\begin{aligned}
&\int_0^\infty |r^{d-1}\Psi(r\lambda^{1/d}) - r^{d-1}\Psi(r\mu^{1/d})|\,dr \\
&\quad \leq \int_0^\infty r^{d-1}[M(r\mu^{1/d}) - M(r\lambda^{1/d})]\,dr \\
&\qquad + \int_0^\infty r^{d-1}[N(r\mu^{1/d}) - N(r\lambda^{1/d})]\,dr \\
&\quad = (M_0 + N_0)(\lambda^{-1} - \mu^{-1}) =: C(\lambda^{-1} - \mu^{-1}).
\end{aligned}$$

Now as in the proof in Example 1.2 of Giné, Mason and Zaitsev [21] choose open balls with centers at $\lambda_k = C/(C - k\varepsilon)$ for $k = 0, \ldots, k_0$ where $k_0$ is the largest integer strictly less than $C/\varepsilon$. This shows (2.24).

Now consider the class of functions on $\mathbf{R}^d$ given by $\mathcal{K} = \{x \mapsto K(\lambda^{1/d}x) : \lambda \geq 1\}$, where $K(\lambda^{1/d}x) = \Psi(\lambda^{1/d}|x|)$. Since by changing to polar coordinates

$$\int_{R^d} |K(\lambda^{1/d}x) - K(\mu^{1/d}x)|\,dx = C_d \int_0^\infty r^{d-1}|\Psi(r\lambda^{1/d}) - \Psi(r\mu^{1/d})|\,dr,$$

where $C_d = d\pi^{d/2}/\Gamma(1 + d/2)$, we see from (2.24) that for some $B_1 > 0$

$$N(\mathcal{K}_\lambda, d_1, \varepsilon) \leq B_1 \varepsilon^{-1}, \qquad 0 < \varepsilon \leq 1. \tag{2.25}$$

From this result along with boundedness of $K$ we readily get that for all $p \geq 1$ there is a $B_p > 0$ such that

$$N(\mathcal{K}_\lambda, d_p, \varepsilon) \leq B_p \varepsilon^{-p}, \qquad 0 < \varepsilon \leq 1. \tag{2.26}$$

It then follows easily that the class of functions $\mathcal{K}_{[a,b]}$ formed from a kernel of the form $K(x) = \Psi(|x|)$ obeying the above conditions satisfies the entropy condition of Theorems 2 and, if $K$ has bounded support or decreases exponentially in a positive power of $|x|$, that of Theorem 3 as well. Some commonly used kernels defined on $\mathbf{R}$ of this form are (1) $K(u) = 1\{u \in [-1/2, 1/2]\}$, (2) $K(u) = \frac{1}{2}\exp(-|u|)$, (3) $K(u) = \frac{1}{\sqrt{2\pi}}\exp(-u^2/2)$, and (4) $K(u) = \frac{3}{4}(1 - x^2)_+$. See, for instance, Devroye [8].

2.2. *Examples.* In this subsection we show how the previous theorems apply to a few instances of estimation of the density of a function of several variables considered in the literature, as well as to the interpoint distance.

EXAMPLE 1 (*Linear combinations*). (Frees [15], Schick and Wefelmeyer [32]) Suppose that

$$g(x_1,\dots,x_m)=\sum_{i=1}^m u_i(x_i),\qquad x_i\in S, \tag{2.27}$$

for measurable functions $u_1,\dots,u_m$ from $S$ to $\mathbf{R}^d$ such that the random variable $u_i(X)$ has a density $f_i$ for each $i=1,\dots,m$. Then (CD) holds with $\overline{f}_i(t,x)=\tilde{f}_i(t-u_i(x))$ and $\tilde{f}_i=f_1*\cdots*f_{i-1}*f_{i+1}*\cdots*f_m$. Thus,

$$\overline{f}(t,x)=\sum_{i=1}^m \tilde{f}_i(t-u_i(x)),\qquad t\in\mathbf{R}^d, x\in S.$$

The process $u_{n,\lambda}$ is then given by

$$\begin{aligned} u_{n,\lambda}(t)=\sqrt{n}\frac{(n-m)!}{n!}\sum_{\mathbf{i}\in I_n^m}\Biggl[&K_{\lambda h_n}\Biggl(t-\sum_{r-1}^m u_r(X_{i_r})\Biggr)\\ &-EK_{\lambda h_n}\Biggl(t-\sum_{r-1}^m u_r(X_r)\Biggr)\Biggr]. \end{aligned}\tag{2.28}$$

We will discuss some conditions under which the central limit theorems in $L_p(\mathbf{R}^d)$, $p\in[1,\infty]$, uniform in $\lambda$, given above apply in this situation.

Let $d=1$. If $f$ is of bounded variation on $\mathbf{R}$ then so is the convolution of $f$ with any density, as is easy to check (in fact, more is true; see, e.g., Schick and Wefelemeyer [33], Lemma 1). Hence $\tilde{f}_j$ is of bounded variation for each $j$ if at least two of the densities $f_i$ are. So, assuming the densities $f_i$ are of bounded variation, then the classes $\mathcal{F}_i=\{\tilde{f}_i(t-\cdot):t\in\mathbf{R}^d\}$ are of VC type by Lemma 22 of Nolan and Pollard [26], as mentioned in connection with (2.12). Also, since the map $(t,x)\mapsto\tilde{f}_i(t-x)$ is jointly measurable, these classes are $Q$-Donsker for every $Q$ (e.g., Dudley [10], Theorem 6.3.1, page 208). Hence, $\mathcal{F}_i$ is $P$-Donsker for each $i=1,\dots,m$. As observed, for example, in the proof of Lemma 8, Schick and Wefelmeyer [34], if $f$ is a function of bounded variation on $\mathbf{R}$, with $f(-\infty+)=0$ and $T_f$ is the total variation of its right continuous modification $f(\cdot+)$, then by Fubini's theorem applied to the product of Lebesgue measure with the measure whose cumulative distribution function is $f(\cdot+)$, $\int |f(x+s)-f(x)|\,dx \le C|s|$ for all $s\in\mathbf{R}$. With $f=\tilde f_i$ this gives

$$E(\overline{f}_i(t,X)-\overline{f}_i(s,X))^2=\int(\tilde f_i(t-u)-\tilde f_i(s-u))^2 f_i(u)\,du \le 2T_{f_i}\|\tilde f_i\|_\infty\|\tilde f_i\|_1|t-s|.$$

Hence, since this holds for every $i=1,\dots,m$, the identity map $(\mathbf{R},|\cdot|)\mapsto(\mathbf{R},\rho)$ is uniformly continuous. So, if the densities $f_i$ are of bounded variation on $\mathbf{R}$, then condition (b) in Theorem 1 is satisfied (note that $f_g$ is of bounded variation, hence bounded).

Let $d>1$. Finding interesting conditions on the densities $f_i$ in order for condition (b) in Theorem 1 in the CLT to be satisfied is a little more cumbersome in this case. Here are some conditions:

(1) $f_i$ is $\alpha$-Hölder continuous for some $0<\alpha\le 1$ and of bounded support, for each $1\le i\le m$. In this case one can easily check the VC property.

(2) $f_i$ is of bounded variation on $\mathbf{R}^d$ in the sense that it is the difference of the ($d$-dimensional) distribution functions of two positive measures (necessarily of the same mass). In this case the Donsker property follows from Dudley [10], Corollary 10.2.8, page 327.

(3) The functions $\tilde f_i$ satisfy condition (2.12). This condition is directly imposed on convolutions because (2.12) may not be inherited by convolution, except in particular cases. Some important examples, like full $d$-variate normal densities, satisfy it.

As a consequence of the above discussion, we have proved the following theorem that improves on Theorem 1 in Schick and Wefelmeyer [32] in that the convergence is uniform in $\lambda\in[a,b]$, the kernel is not necessarily a convolution and the window sizes $h_n$ are allowed to decrease at a smaller rate. Their result includes the bias part whereas ours does not, and we discuss this immediately below.

THEOREM 4. *Let $g$ be defined by (2.27) and let $K$ satisfy condition (2.12). If $d=1$ assume that the densities $f_i$ of $u_i(X)$ are of bounded variation and if $d>1$ assume that $f_i$ satisfy* (1) *or* (2) *or that $\tilde f_i$ satisfy* (3) *in the previous paragraph. Let $h_n\to 0$ and $nh_n/(\log h_n^{-1})^2\to\infty$. Then the processes $u_{n,\lambda}$ defined by (2.28) converge weakly in $L_\infty(\mathbf{R}^d)$ uniformly in $a\le\lambda\le b$ to the centered Gaussian process with the same covariance as $\overline{f}(\cdot,x)=\sum_{i=1}^m \tilde f_i(\cdot-u_i(x))$.*

Next we comment on the bias part in this theorem. What we do is standard and can also be done for the rest of the results in this article, but we will refrain from doing so. Suppose that the densities $f_i$ are in $C^k(\mathbf{R}^d)$ and that their partial derivatives of order $k$ or smaller are all bounded. Then the same is true for $f_g = f_1 * \cdots * f_m$ (this follows, e.g., from Proposition 8.10 in Folland [14] and Young's inequalities) and therefore $f_g$ admits a Taylor development of the form

$$f_g(t+\delta) - f_g(t) = \sum_{r=1}^{k} \frac{1}{r!} \sum_{\substack{s_1+\cdots+s_d=r \\ 0\le s_i\le r}} \frac{\partial^r f_g}{\partial^{s_1} x_1 \cdots \partial^{s_d} x_d}(t)\delta_1^{s_1}\cdots\delta_d^{s_d} + H(t,k,\delta)|\delta|^k,$$

where $H$ is uniformly bounded by a constant times the common bound for the $k$th partial derivatives of $f_g$ (actually, by continuity, $H \to 0$ as $\delta \to 0$). Suppose moreover that the kernel $K$ satisfies the conditions

$$\int |t|^k |K(t)|\,dt < \infty \quad \text{and} \quad \int t_1^{s_1}\cdots t_d^{s_d} K(t)\,dt = 0 \qquad \text{for } s_i \ge 0, \sum_{i=1}^{d} s_i \le k.$$

Then, integration after change of variables and use of Taylor's formula for $f_g(t+(\lambda h_n)^{1/d}u) - f_g(t)$ give

$$\sup_{\lambda\in[a,b]} \sqrt{n}|EK_{\lambda h_n}(t - g(X_1,\ldots,X_m)) - f_g(t)| = O(\sqrt{n}h_n^{k/d}).$$

Therefore, under these extra conditions on $f_i$ and $K$, if $nh_n^{2k/d} \to 0$, the conclusion of the previous theorem can be modified to convergence in law in $\ell^\infty(\mathbf{R}^d \times [a,b])$ of

$$\sqrt{n}\left\{\frac{1}{|I_n^m|}\sum_{I_n^m} K_{\lambda h_n}\left(t - \sum_{r=1}^{m} u_r(X_{i_r})\right) - f_g(t)\right\}$$

to the same centered Gaussian process. If the partial derivatives up to order $k$ of $f_i$ are bounded and in $L_1$, then $f_g \in C^{mk}(\mathbf{R}^d)$ (e.g., by iteration in Proposition 8.10, Folland [14]), and the previous discussion applies with $k$ replaced by $mk$.

*Simultaneous estimation of convolutions.* Here is an immediate application of Theorem 4 to the simultaneous estimation of the densities of convolutions. Let $f$ be a density of bounded variation on $\mathbf{R}$. We are interested in estimating the convolutions $f^{*2} = f * f, \ldots, f^{*N}$, for $N \ge 2$. For $m = 2, \ldots, N$, $t_m \in \mathbf{R}$ and $a_m \le \lambda_m \le b_m$ introduce the estimators

$$\widehat{f_n^{*m}}(t_m,\lambda_m) = \frac{1}{|I_n^m|\lambda_m h_n} \sum_{I_n^m} K\left(\frac{t_m - \sum_{r=1}^m X_{i_r}}{\lambda_m h_n}\right)$$

and set

$$u_{n,\lambda_m}^{(m)}(t_m) = \sqrt{n}\{\widehat{f_n^{*m}}(t_m,\lambda_m) - E\widehat{f_n^{*m}}(t_m,\lambda_m)\}$$

or, under extra conditions as in the previous remark, replace $E\widehat{f_n^{*m}}(t_m,\lambda_m)$ by $f^{*m}(t_m,\lambda_m)$. [These estimators correspond to the functions $g=g_m$ defined by $g_m(x_1,\dots,x_m)=x_1+\cdots+x_m$.] A direct application of Theorem 4 in combination with the observations (2.6)–(2.8) gives that the vector-valued processes $(u^{(2)}_{n,\lambda_2}(t_2),\dots,u^{(N)}_{n,\lambda_N}(t_N))$ converge weakly in $L_\infty(\mathbf{R})\times\cdots\times L_\infty(\mathbf{R})$ uniformly in $a_m\le\lambda_m\le b_m$, $m=2,\dots,N$, to the centered vector-valued Gaussian process defined on $\mathbf{R}\times\cdots\times\mathbf{R}$, with the same covariance/cross covariance matrix as

$$(2f(t_2-X),\dots,Nf^{*(N-1)}(t_N-X)).$$

We note here that Schick and Wefelmeyer [32] use a variation of the Frees [15] local $U$-statistic estimator of convolutions of densities. Their estimator is based on convolving kernel density estimators. Here is how their approach works in the case of estimating the density of $X_1+X_2$, where $X_1$ and $X_2$ are i.i.d. real valued with density $f$. Consider the kernel density estimator

$$\widehat{f}_n(x)=(nh_n)^{-1}\sum_{i=1}^n k(h_n^{-1}(x-X_i)),$$

where $k$ is of bounded variation on $\mathbf{R}$. Their estimator is $\widehat{f}_n*\widehat{f}_n$ and can be expressed, with $K=k*k$, as

$$\begin{aligned}\widehat{f}_n*\widehat{f}_n(x)&=\frac{1}{n^2h_n}\sum_{i=1}^n\sum_{j=1}^n K\left(\frac{x-X_i-X_j}{h_n}\right)\\&=\frac{n-1}{n}\overline{U}_{n,h_n}(x)+\frac{1}{n^2h_n}\sum_{i=1}^n K\left(\frac{x-2X_i}{h_n}\right),\end{aligned}$$

where $\overline{U}_{n,h_n}(x)$ is the Frees type local $U$-statistic estimator of $f*f(x)$ defined by

$$\overline{U}_{n,h_n}(x)=(h_nn(n-1))^{-1}\sum_{1\le i\ne j\le n}K(h_n^{-1}(x-X_i-X_j)).$$

The second term in the above expression for $\widehat{f}_n*\widehat{f}_n(x)$ is asymptotically negligible and Theorem 4 applies to $\overline{U}_{nh_n}$. This remark applies as well to simultaneous estimation of convolutions of densities.

We can complement the above results for $L_\infty(\mathbf{R}^d)$ with limit theorems for the $L_p$ distance. We now do this for the cases $p=1$ and $p=2$.

THEOREM 5. *Assume $K$ satisfies condition* (a) *in Theorem 2 and let $g(x_1,\dots,x_m)=\sum_{i=1}^m u_i(x_i)$, $u_i:S\mapsto\mathbf{R}^d$ be measurable, where $u_i(X)$ has* *density $f_i$. Assume $f_i$ is in $L_2(\mathbf{R}^d)$. Assume also $h_n \to 0$ and $nh_n \to \infty$. Then*

$$\frac{\sqrt{n}}{|I_n^m|}\sum_{I_n^m}\left\{K_{\lambda h_n}\left(t-\sum_{r=1}^m u_r(X_{i_r})\right)-EK_{\lambda h_n}\left(t-\sum_{r=1}^m u_r(X_r)\right)\right\}$$

*converges in law in $L_2(\mathbf{R}^d)$ uniformly in $\lambda \in [a,b]$ to the centered Gaussian process with the same covariance as $\sum_{i=1}^m \tilde{f}_i(\cdot - u_i(X_i))$.*

This theorem follows because by Young's inequality (e.g., Folland [14], page 240), $f_i \in L_2$ for $i = 1,\dots,m$ implies $E\int \overline{f}(t,X)^2\,dt < \infty$. As in the previous theorem, a little more smoothness on $f_i$ and higher-order kernels allow for elimination of the bias.

We can also recover the Schick and Wefelmeyer [32], Theorem 2 on the CLT for the $L_1$ norm, with weaker assumptions (note that condition (2.22) is vacuous if $[a,b]$ reduces to a single point).

If two densities belong to $L_1(\mu_s)$ for some $s > 0$, then so does their convolution as can be seen by direct computation using the trivial observation that $(1+|u+v|) \le (1+|u|)(1+|v|)$, and it is also routine to check that if any finite number of densities and their squares belong to $L_1(\mu_s)$, so does the square of their convolution. Moreover, by Lemma 1 to be proved below, if $f$ and $k$ are two nonnegative functions in $L_1(\mu_s)$ for some $s > d$, then $\int (f*k)^{1/2}(t)\,dt < \infty$. Hence Theorem 3 gives:

THEOREM 6. *Let $g(x_1,\dots,x_m) = \sum_{i=1}^m u_i(x_i)$, $u_i : S \mapsto \mathbf{R}^d$ be measurable, with $u_i(X)$ having density $f_i$, for $i = 1,\dots,m$, and let $K$ be a kernel on $\mathbf{R}^d$. Assume that for some $s > d$, $K^2$, $f_i$ and $f_i^2$ are in $L_1(\mu_s)$, $i = 1,\dots,m$, and that $K$ satisfies condition* (a) *in Theorem 3 for this $s$. Assume also $h_n \to 0$ and $nh_n \to \infty$. Then*

$$\frac{\sqrt{n}}{|I_n^m|}\sum_{I_n^m}\left\{K_{\lambda h_n}\left(t-\sum_{r=1}^m u_r(X_{i_r})\right)-EK_{\lambda h_n}\left(t-\sum_{r=1}^m u_r(X_r)\right)\right\}$$

*converges in law in $L_1(\mathbf{R}^d)$ uniformly in $\lambda \in [a,b]$) to the centered Gaussian process with the same covariance as $\sum_{i=1}^m \tilde{f}_i(\cdot - u_i(X_i))$.*

In $\mathbf{R}$, if the densities $f_i$ are bounded (they do not need to be) then the condition imposed on $f_i$ is simply that $\int |x|^{1+\delta} f(x)\,dx < \infty$ for some $\delta > 0$.

EXAMPLE 2 (*Distribution of sample distances*). Frees [15] considers estimating the density of the interpoint functional $g(X_1,X_2) = |X_1 - X_2|$ in two dimensions, where $X_i$ are i.i.d. with a density $f$ which is bounded and of bounded support. We further assume that $f$ is $\alpha$-Hölder continuous for some $0<\alpha\le 1$, that is, $|f(u)-f(v)|\le C|u-v|^\alpha$ for all $u,v\in\mathbf{R}^2$. Then, it is easy to see that $f_g$, the density of $g$,

$$f_g(t)=\int_0^{2\pi}\int_{\mathbf{R}^2} f(x_1+t\cos\theta,x_2+t\sin\theta)f(x_1,x_2)t\,dx_1\,dx_2\,d\theta,$$

is bounded and has bounded support. Moreover, for $x=(x_1,x_2)$ and with $R$ the radius of a ball around zero containing the support of $f$, we have, for the conditional densities $\overline{f}_i$ of $|X_1-X_2|$,

$$\overline{f}_1(t,x)=\overline{f}_2(t,x)=\int_0^{2\pi} f(x_1+t\cos\theta,x_2+t\sin\theta)t\,d\theta,\qquad 0\le t\le 2R,$$

and $\overline{f}_i(t,x)=0$ for larger values of $t$. Then $\overline{f}=\overline{f}_1+\overline{f}_2$ satisfies

$$|\overline{f}(t,\cdot)-\overline{f}_i(s,\cdot)|\le 8\pi RC|t-s|^\alpha+4\pi\|f\|_\infty|t-s|,$$

and therefore, with $\mathcal{F}=\{\overline{f}(t,\cdot):|t|\le 2R\}$ and any probability measure $Q$, $N(\mathcal{F},L_2(Q),\varepsilon)\le C/\varepsilon^{1/\alpha}$ which, since the class $\mathcal{F}$ is image admissible Suslin by joint measurability of $\overline{f}$, implies that the class $\mathcal{F}$ is $P$-Donsker. So, if we take a kernel $K$ satisfying (2.12), and $h_n\to 0$ with $nh_n/(\log h_n^{-1})^2\to\infty$, we get that the processes defined for $|t|\le R$, $\lambda\in[a,b]$,

$$\frac{n^{1/2}}{n(n-1)}\sum_{1\le i\ne j\le n}\{K_{\lambda h_n}(t-|X_i-X_j|)-EK_{\lambda h_n}(t-|X_i-X_j|)\},$$

converge in law uniformly in $t$ and $\lambda$ to the centered Gaussian process with the same covariance as $2f_1(t,X)$. Moreover, the comments following Theorem 4 regarding replacement of $(\lambda h_n)^{-1}EK((\lambda h_n)^{-1}(t-|X_i-X_j|))$ by $f_g(t)$ apply here as well.

In dimension 1,

$$f_1(t,x)=f_2(t,x)=f(t+x)+f(-t+x),\qquad t\ge 0,$$

and a sufficient condition for $\mathcal{F}$ to be $P$-Donsker is that $f$ be of bounded variation on $\mathbf{R}$. Then, $K$ of bounded variation on $\mathbf{R}$ and $nh_n/(\log h_n^{-1})^2\to\infty$ ensure the same CLT as above. Since

$$f_g(x)=\int f(y+x)f(y)\,dy+\int f(y-x)f(y)\,dy$$

and

$$f_g(0)=2\int f^2(y)\,dy,$$

as observed by Frees [15], we get a $\sqrt{n}$ consistent estimator of $\int f^2(x)\,dx$ when the extra smoothness conditions to make the bias tend to zero hold. See also Bickel and Ritov [4] or Giné and Mason [19] and references therein for other $\sqrt{n}$ consistent estimators of $\int f^2(x)\,dx$.

EXAMPLE 3 (*Local interpoint distance processes*). (Jammalamadaka and Janson [23], Eastwood and Horváth [12].) In Corollary 1 we consider the processes

$$\frac{\sqrt{n}}{h_n \operatorname{Vol}(D)} \sum_{I_n^m} [I\{g(X_{i_1}, \ldots, X_{i_m}) \in (\lambda h_n)^{1/d} D\} \\ - \Pr\{g(X_1, \ldots, X_m) \in (\lambda h_n)^{1/d} D\}]$$

for $D$ star shaped about 0. Of particular interest in the literature is the case corresponding to $m = 2$, $S = \mathbf{R}^d$, $g(X_1, X_2) = X_2 - X_1$ and $D$ the (open or closed) unit ball about zero for some norm in $\mathbf{R}^d$. In this case the densities $\overline{f}_i(t, x)$, $i = 1, 2$, are, respectively, $f(x - t)$ and $f(x + t)$, where $f$ is the density of $X$. So, for the local asymptotic equicontinuity condition (2.13) to hold we only need that the class $\mathcal{F} = \{f(\cdot + t) : |t| < \delta\}$ be $P$-Donsker for some $\delta > 0$. If $f$ is Hölder continuous of order $\alpha \in (0, 1]$, and $Q$ is any probability measure on $\mathbf{R}^d$, then

$$\int (f(x + t) - f(x + s))^2 \, dQ(x) \le C|t - s|^{2\alpha}, \qquad s, t \in \mathbf{R}^d.$$

It follows as in the previous example that $\mathcal{F}$ is VC and measurable, hence $Q$-Donsker for every probability measure $Q$. Thus, Corollary 1 implies the following slight strengthening and generalization of Theorem 1.1 in Eastwood and Horváth [12].

THEOREM 7. *If $D$ is a bounded measurable subset of $\mathbf{R}^d$ star-shaped about zero and with $\operatorname{Vol}(D) \neq 0$, $X_i$ are i.i.d. random vectors in $\mathbf{R}^d$ with a density $f$ which is $\alpha$-Hölder continuous for some $\alpha \in (0, 1]$, and $h_n \to 0$, $nh_n \to \infty$, then the processes*

$$\frac{1}{n^{3/2} h_n \operatorname{Vol}(D)} \sum_{(i_1, i_2) \in I_n^2} \left[ I\left\{ \frac{X_{i_1} - X_{i_2}}{(\lambda h_n)^{1/d}} \in D \right\} - \Pr\left\{ \frac{X_1 - X_2}{(\lambda h_n)^{1/d}} \in D \right\} \right], \qquad 0 \le \lambda \le 1,$$

*converge in law in $\ell^\infty[0, 1]$ to the process $\lambda \sigma Z$, $0 \le \lambda \le 1$, where $Z$ is $N(0, 1)$ and $\sigma^2 = 4[\int f^3(x) \, dx - (\int f^2(x) \, dx)^2]$.*

## 3. Proofs.

In the sequel it will be helpful to introduce the following notation and facts. For a kernel $L$ of $k \ge 1$ variables we set

$$U_n^{(k)}(L) = \frac{(n-k)!}{n!} \sum_{\mathbf{i} \in I_n^k} L(X_{i_1}, \ldots, X_{i_k}) \tag{3.1}$$

[so, $U_n(t, \lambda) = U_n^{(m)}(K_{\lambda h_n}(t - g(\cdot, \ldots, \cdot)))$ even if $g$ is not symmetric in its entries]. [When $L$ is a constant function we define $U_n^{(0)}(L) = L$.] Assume now that $L$ is a function of $m \geq 1$ variables, symmetric in its entries. Then, for $1 \leq k \leq m$, the Hoeffding projections with respect to $P$ are defined as

$$\pi_k L(x_1, \ldots, x_k) = (\delta_{x_1} - P) \times \cdots \times (\delta_{x_k} - P) \times P^{m-k}(L) \tag{3.2}$$

and $\pi_0 L = EL(X_1, \ldots, X_m)$. Then, the *Hoeffding decomposition* states the following, which is easy to check:

$$U_n^{(m)}(L) - EL = \sum_{k=1}^{m} \binom{m}{k} U_n^{(k)}(\pi_k L). \tag{3.3}$$

For $L \in L_2(P^m)$ this is an orthogonal decomposition and $E(\pi_k L | X_2, \ldots, X_k.) = 0$ for $k \geq 1$; that is, the kernels $\pi_k L$ are canonical for $P$ (or completely degenerate, or completely centered). Also, $\pi_k$, $k \geq 1$, are nested projections, that is, $\pi_k \circ \pi_\ell = \pi_k$ if $k \leq \ell$, and $E(\pi_k L)^2 \leq E(L - EL)^2 \leq EL^2$.

The function $K_h(t - g(X_1, \ldots, X_m))$ is not necessarily symmetric in its entries, but we can symmetrize it as

$$\overline{K}_h(t, x_1, \ldots, x_m) := \frac{1}{m!} \sum_{\sigma \in I_m^m} K_h(t - g(x_{\sigma_1}, \ldots, x_{\sigma_m})). \tag{3.4}$$

Then, clearly, for each $t \in \mathbf{R}^d$,

$$\begin{aligned} &U_n(t, \lambda) - EK_{\lambda h_n}(t - g(X_1, \ldots, X_m)) \\ &\quad = U_n^{(m)}(\overline{K}_{\lambda h_n}(t, \cdot, \ldots, \cdot)) - E\overline{K}_{\lambda h_n}(t, X_1, \ldots, X_m). \end{aligned}$$

Moreover, by applying (3.3) to $u_{n,\lambda}(t)$ we get

$$u_{n,\lambda}(t) = \sqrt{n} \sum_{k=1}^{m} \binom{m}{k} U_n^{(k)}(\pi_k \overline{K}_{\lambda h_n}(t, \cdot, \ldots, \cdot)). \tag{3.5}$$

3.1. *A general proposition for the CLT.* Let us consider the first term in the expansion (3.5). Note that, by definition of $\overline{f}$ and $f_g$,

$$\begin{aligned} m\pi_1 \overline{K}_{\lambda h_n}(t, x) &= \sum_{i=1}^{m} EK_{\lambda h_n}(t - g(X_1, \ldots, X_{i-1}, x, X_{i+1}, \ldots, X_m)) \\ &\quad - mEK_{\lambda h_n}(t - g(X_1, \ldots, X_m)) \\ &= \int K_{\lambda h_n}(t - u)(\overline{f}(u, x) - E\overline{f}(u, X))\, du \\ &= \int (\overline{f}(t - u, x) - E\overline{f}(t - u, X)) K_{\lambda h_n}(u)\, du. \end{aligned}$$

Hence,

$$\sqrt{n} m U_n^{(1)}(\pi_1 \overline{K}_{\lambda h_n}(t, \cdot))$$

$$= \frac{1}{\sqrt{n}} \sum_{i=1}^{n} \int (\overline{f}(t-u, X_i) - E\overline{f}(t-u, X)) K_{\lambda h_n}(u)\, du \tag{3.6}$$

$$= \int \overline{\nu}_n(t-u) K_{\lambda h_n}(u)\, du = (\overline{\nu}_n * K_{\lambda h_n})(t),$$

which is a generalized version of the smoothed empirical process that has been recently investigated by several authors (e.g., Rost [31] and references therein). We shall see that it controls the asymptotic behavior of the process $u_{n,\lambda}$. In the next proposition we show it has the same asymptotic behavior as the empirical process over the class of functions $\{\overline{f}(t,\cdot), t \in \mathbf{R}^d\}$.

PROPOSITION 1. *Let $1 \le p \le \infty$. Assume:*

(i) *condition (CD$p$) holds;*

(ii) $\lim_{\delta\to 0} \limsup_{n\to\infty} \Pr^*\{\sup_{|u|\le\delta} \|\overline{\nu}_n(\cdot - u) - \overline{\nu}_n(\cdot)\|_p > \varepsilon\} = 0$ *for all* $\varepsilon > 0$*;*

(iii) *the sequence $\|\overline{\nu}_n\|_p^*$, $n \in \mathbf{N}$, is stochastically bounded.*

*Then, whenever $h_n \to 0$ we have*

$$\lim_{n\to\infty} \sup_{a \le \lambda \le b} \|\sqrt{n} m U_n^{(1)}(\pi_1 \overline{K}_{\lambda h_n}(\cdot,\cdot)) - \overline{\nu}_n(\cdot)\|_p = 0 \qquad \textit{in } \mathrm{pr}^*. \tag{3.7}$$

PROOF. Let $w_\delta(\overline{\nu}_n) = \sup_{|u|<\delta} \|\overline{\nu}_n(\cdot - u) - \overline{\nu}(\cdot)\|_p$, $\delta > 0$, denote the $L_p(\mathbf{R}^d)$ modulus of "continuity" of $\overline{\nu}_n$, which is defined because of (i). Then it follows by Fubini in the case $p = \infty$ and by Minkowski's inequality for integrals (e.g., Folland [14], page 194) in the case $1 \le p < \infty$ that

$$\|\overline{\nu}_n * K_{\lambda h_n} - \overline{\nu}_n\|_p \le w_\delta(\overline{\nu}_n)\|K_{\lambda h_n}\|_1 + 2\|\overline{\nu}_n\|_p \int_{|u|>\delta} |K_{\lambda h_n}(u)|\, du$$

$$\le w_\delta(\overline{\nu}_n)\|K\|_1 + 2\|\overline{\nu}_n\|_p \int_{|u|>\delta/(\lambda h_n)^{1/d}} |K(u)|\, du.$$

Now, the result follows from this, $K \in L_1(\mathbf{R}^d)$ and (ii) and (iii), in view of (3.6). □

COROLLARY 2. *Let $1 \le p \le \infty$. Assume (CD$p$) and hypothesis* (b) *in Theorem 1 for $p = \infty$, in Theorem 2 for $2 \le p < \infty$ and in Theorem 3 for $1 \le p < 2$. Then the processes $\overline{\nu}_n$ converge weakly in $L_p(\mathbf{R}^d)$ to the centered Gaussian process $G$ with the covariance of $\overline{f}(\cdot, X)$. If moreover $h_n \to 0$, then the limit (3.7) holds.*

PROOF. (a) Case $p = \infty$. In this case, the hypothesis of the class $\mathcal{F}$ being $P$-Donsker is just another way of saying that $\overline{\nu}_n$ converges in law to $G$ in $L_\infty(\mathbf{R}^d)$, and this obviously implies (iii); finally, (ii) holds by the uniform continuity of the identity map $(\mathbf{R}^d, |\cdot|) \mapsto (\mathbf{R}^d, \rho)$ together with the usual asymptotic equicontinuity condition (e.g., Dudley [10], pages 117–118).

(b) Case $1 \le p < \infty$. As mentioned in Section 2, the hypotheses (b) in Theorems 2 and 3 are precisely the necessary and sufficient conditions for the process $\overline{f}(\cdot, X)$ to satisfy the CLT in $L_p$, that is, for $\overline{\nu}_n$ to converge in law to $G$ in $L_p$. Moreover, (iii) is a direct consequence of this convergence. Finally, condition (ii) holds by the uniform tightness implied by weak convergence together with the Fréchet–Kolmogorov characterization of compact sets of $L_p(\mathbf{R}^d)$ (see, e.g., Dunford and Schwartz [11], Theorem IV.8.21, page 301). □

In view of (3.5) and Corollary 2, to complete our CLT program for $u_{n,\lambda}$, that is, the proofs of Theorems 1, 2 and 3, it only remains to show that the hypotheses in the statements of these theorems also imply

$$\sup_{a\le\lambda\le b} \|\sqrt{n} U_n^{(k)}(\pi_k \overline{K}_{\lambda h_n})\|_p \to 0 \qquad \text{in pr}, \quad k = 2, \dots, m. \tag{3.8}$$

Proving this constitutes the main part of our proofs, and requires new inequalities for $U$-processes that we develop in the next subsections.

3.2. *Inequalities for $U$-processes.* In the next two subsections we collect the inequalities we need to prove the limits (3.8). First we consider the case of $U$-processes indexed by VC classes of functions and obtain a moment inequality that generalizes the scope of that of Einmahl and Mason [13] and Giné and Koltchinskii [17] for empirical processes. Next, we consider $U$-statistics taking values on separable type 2 Banach spaces such as $L_p$, $p \ge 2$, and derive exponential inequalities for them. The inequality that we get in this situation is particularly neat. Then, based on the method recently used by Giné, Latała and Zinn [18] to prove inequalities for $U$-statistics, we derive both moment and exponential inequalities for other Banach spaces, such as $L_p$, $1 \le p < 2$. These are less clean than in the type 2 case, but are still usable. Neither of our exponential inequalities captures the Gaussian tail behavior that the statistic should have for small values of $x$; nevertheless, their application yields very strong results in the situations encountered in this article.

3.3. *$U$-processes indexed by VC classes.* In this subsection we consider classes of measurable functions $\mathcal{F}$ defined on $(S^m, \mathcal{S}^m)$ taking values in $[-1, 1]$, and we assume that $0 \in \mathcal{F}$. The object is to obtain a bound for $E\|U_n^{(k)}(\pi_k f)\|_{\mathcal{F}}$ where $\mathcal{F}$ is of VC-type, and where we use the notation $\|\Psi(f)\|_{\mathcal{F}} = \sup_{f\in\mathcal{F}} |\Psi(f)|$ for any functional $\Psi$ defined on the class $\mathcal{F}$. This bound will require measurability on the class $\mathcal{F}$ described in de la Peña and Giné [7], page 138: the class $\mathcal{F}$ should be measurable in the sense that for every $k = 1, \dots, m$ and every choice of $a_{i_1,\dots,i_k} \in \{-1, 1\}$, the mapping $f \mapsto \sum_{I_n^k} a_{i_1,\dots,i_k} P^{m-k} f(X_{i_1}, \dots, X_{i_k})$ is measurable for the completion of $\mathcal{S}^n$. This holds, for instance, if (a) if there exists $\mathcal{F}_0$ countable such that this sup equals the sup over $\mathcal{F}_0$, or (b) if the $\sigma$-algebra $\mathcal{S}$ is countably generated and contains the singletons, and the class $\mathcal{F}$ is image admissible Suslin, for instance, if it is parametrized by a complete separable metric space $T$ in such a way that the evaluation map $(t, x_1, \dots, x_m) \mapsto f_t(x_1, \dots, x_m)$ is jointly measurable (Dudley [10], Section 5.3; van der Vaart and Wellner [37], Section 2.3.1; Pollard [30], page 196). These conditions allow us to randomize by independent random signs and use Fubini. If either of these two conditions is satisfied, we say that the class $\mathcal{F}$ is measurable.

The following moment inequality will be instrumental in finishing the proof of Theorem 1. The proof of this inequality has several points in common with the proofs of similar inequalities for $m = 1$ in Einmahl and Mason [13] and in Giné and Koltchinskii [17]; however the present proof does not rely on the square root trick or on the contraction principle for Rademacher processes.

THEOREM 8. *Let $\mathcal{F}$ be a measurable collection of functions $S^m \mapsto \mathbf{R}$ symmetric in their entries with an envelope function $F$ and let $P$ be any probability measure on $(S, \mathcal{S})$ (with $X_i$ i.i.d. $P$). Assume $F$ is bounded by $M > 0$ and $\mathcal{F}$ is VC with respect to $F$ with characteristics $A$ and $v$, as in (2.9). Then for every $m \in \mathbf{N}$, $A \geq e^m$, $v \geq 1$, there exist constants $C_1 := C_1(m, A, v, M)$ and $C_2 = C_2(m, A, v, M)$ such that*

$$n^{k/2} E\|U_n^{(k)}(\pi_k f)\|_{\mathcal{F}} \leq C_1 \sigma \left(\log \frac{A\|F\|_{L_2(P^m)}}{\sigma}\right)^{k/2}, \qquad k = 0, 1, \dots, m, \tag{3.9}$$

*assuming*

$$n\sigma^2 \geq C_2 \log\left(\frac{2\|F\|_{L_2(P^m)}}{\sigma}\right), \tag{3.10}$$

*where $\sigma^2$ is any number satisfying*

$$\|P^m f^2\|_{\mathcal{F}} \leq \sigma^2 \leq P^m F^2. \tag{3.11}$$

PROOF. Without loss of generality we assume $F \leq M = 1$. The theorem is true for all $m$ and $k = 0$ by Hölder's inequality, since $U_n^{(0)}(\pi_0 f) = Pf$, so we assume the statement to be true for all $m \in \mathbf{N}$ and $k - 1$, for some $k \geq 1$, and prove it for $k$ (and for all $m \in \mathbf{N}$ ). We shall omit symbols when no confusion is possible, so, for instance, we write $\|\cdot\|$ for $\|\cdot\|_{\mathcal{F}}$. We shall not keep track of constants; in particular, constants that depend on a subset of $m, A, v, M$, will generically be denoted by $C$ (so, the value of $C$ may change from line to line). By Theorem 3.5.3 in de la Peña and Giné [7], the definition of $\pi_k$ and Jensen's inequality,

$$E\left\|\sum_{I_n^k}(\pi_k f)(X_{i_1},\dots,X_{i_k})\right\| \le CE\left\|\sum_{I_n^k}\varepsilon_{i_1}\cdots\varepsilon_{i_k}(\pi_k f)(X_{i_1},\dots,X_{i_k})\right\|$$

$$\le CE\left\|\sum_{I_n^k}\varepsilon_{i_1}\cdots\varepsilon_{i_k}(P^{m-k} f)(X_{i_1},\dots,X_{i_k})\right\|.$$

Here $X_i, \varepsilon_j$ are all independent, the $X$'s have law $P$ and the $\varepsilon$'s are random signs (Rademacher variables). Since for any probability measure $Q$ on $S^k$,

$$Q(P^{m-k}(f-g))^2 \le Q \times P^{m-k}(f-g)^2,$$

it follows from the VC property of $\mathcal{F}$ that

$$N(P^{m-k}\mathcal{F}, L_2(Q), \tau) \le \left(\frac{A\|\sqrt{P^{m-k}F^2}\|_{L_2(Q)}}{\tau}\right)^v,$$

$$0 < \tau \le 2\|\sqrt{P^{m-k}F^2}\|_{L_2(Q)},$$

that is, $P^{m-k}\mathcal{F}$ is VC-type with characteristics $A$ and $v$ and envelope $\sqrt{P^{m-k}F^2}$. This gives, by the entropy integral for Rademacher chaos of order $k$ (de la Peña and Giné [7], Corollary 5.1.8, upon noting that exponential Orlicz norms dominate $L_p$ norms up to constants),

$$|I_n^k|^{-1/2}E_\varepsilon\left\|\sum_{I_n^k}\varepsilon_{i_1}\cdots\varepsilon_{i_k}(P^{m-k} f)(X_{i_1},\dots,X_{i_k})\right\|$$

$$\le C\int_0^{\sqrt{\|U_n^{(k)}((P^{m-k}f)^2)\|}}\left(\log\frac{A\sqrt{U_n^{(k)}(P^{m-k}F^2)}}{\tau}\right)^{k/2}d\tau,$$

where $E_\varepsilon$ denotes expectation with respect to the Rademacher variables only. [To apply Corollary 5.1.8 exactly to our process we need that, for $X_1,\dots,X_n$ fixed, the Rademacher chaos process

$$f \mapsto \sum_{I_n^k}\varepsilon_{i_1}\cdots\varepsilon_{i_k}(P^{m-k} f)(X_{i_1},\dots,X_{i_k}), \qquad f \in \mathcal{F},$$

be separable, and this follows by separability of the unit cube of $\mathbf{R}^{|I_n^k|}$ for the Euclidean norm.] Then, by Fubini,

$$|I_n^k|^{1/2}E\|U_n^{(k)}(\pi_k f)\| \le CB, \tag{3.12}$$

where

$$B = E\int_0^{\sqrt{\|U_n^{(k)}((P^{m-k}f)^2)\|}}\left(\log\frac{A\sqrt{U_n^{(k)}(P^{m-k}F^2)}}{\tau}\right)^{k/2}d\tau. \tag{3.13}$$

Decompose the integral $B$ into two parts,

$$(I) := E\left[\int_0^{\sqrt{\|U_n^{(k)}((P^{m-k}f)^2)\|}} \left(\log \frac{A\sqrt{U_n^{(k)}(P^{m-k}F^2)}}{\tau}\right)^{k/2} d\tau\, I_n\right],$$

where

$$I_n = I\{U_n^{(k)}(P^{m-k}F^2) > 4P^m F^2\}$$

and

$$(II) := E\left[\int_0^{\sqrt{\|U_n^{(k)}((P^{m-k}f)^2)\|}} \left(\log \frac{A\sqrt{U_n^{(k)}(P^{m-k}F^2)}}{\tau}\right)^{k/2} d\tau\, I_n^c\right],$$

so that $B = (I) + (II)$. The $(I)$ term is handled by the Arcones [3] exponential inequality (de la Peña and Giné [7], Theorem 4.1.13), which gives

$$\begin{aligned}(3.14)\qquad &\Pr\{U_n^{(k)}(P^{m-k}F^2) > 4P^m F^2\}\\ &\quad\le 4\exp\left\{-\frac{9n(P^m F^2)^2}{2k^2 P^m F^4 + cP^m F^2}\right\} \le 4\exp\left\{-\frac{9nP^m F^2}{2k^2 + c}\right\}\end{aligned}$$

for a constant $c$ that depends only on $k$. In the last inequality we have used $F \le 1$. Since, by change of variables,

$$\begin{aligned}&\int_0^{\sqrt{\|U_n^{(k)}((P^{m-k}f)^2)\|}} \left(\log \frac{A\sqrt{U_n^{(k)}(P^{m-k}F^2)}}{\tau}\right)^{k/2} d\tau\\ &\qquad \le A\int_0^1 (\log u^{-1})^{k/2}\, du \sqrt{U_n^{(k)}(P^{m-k}F^2)},\end{aligned}$$

and $P^k(U_n^{(k)}(P^{m-k}F^2)) \le P^m F^2$, it follows by Hölder's inequality and (3.14) that

$$(I) \le C\|F\|_{L_2(P^m)} \exp\left\{-\frac{9n\|F\|^2_{L_2(P^m)}}{2(2k^2 + c_1)}\right\} \le \frac{D}{\sqrt{n}}. \tag{3.15}$$

As for $(II)$, we note

$$(II) \le E\left[\int_0^{\sqrt{\|U_n^{(k)}((P^{m-k}f)^2)\|}} \left(\log \frac{2A\|F\|_{L_2(P^m)}}{\tau}\right)^{k/2} d\tau\, I_n^c\right]. \tag{3.16}$$

Now by integration we see that for any $0 < c < C$

$$\int_0^c (\log(C/x))^{k/2}\left[1 - \frac{k}{2}(\log(C/x))^{-1}\right] dx = c(\log(C/c))^{k/2},$$

which when $(\log(C/c))^{-1} \le k^{-1}$ gives the inequality

$$\int_0^c (\log(C/x))^{k/2}\, dx \le 2c(\log(C/c))^{k/2}. \tag{3.17}$$

Thus since on $I_n^c$

$$2A\|F\|_{L_2(P^m)}/\sqrt{\|U_n^{(k)}((P^{m-k}f)^2)\|} \geq A \geq e^m > e^k, \tag{3.18}$$

we get from (3.17), (3.18) and (3.16) that

$$(II) \leq 2E\Bigg[\sqrt{\|U_n^{(k)}((P^{m-k}f)^2)\|}\Bigg(\log\frac{2A\|F\|_{L_2(P^m)}}{\sqrt{\|U_n^{(k)}((P^{m-k}f)^2)\|}}\Bigg)^{k/2}\Bigg].$$

Since the function $\sqrt{x}(-\log x)^{k/2}$, $0 < x < 1$, is concave on $(0, e^{-k}]$ and $A \geq e^m > e^k$, this last bound is by Jensen's inequality

$$\leq 2\sqrt{E\|U_n^{(k)}((P^{m-k}f)^2)\|}\Bigg(\log\frac{2A\|F\|_{L_2(P^m)}}{\sqrt{E\|U_n^{(k)}((P^{m-k}f)^2)\|}}\Bigg)^{k/2}. \tag{3.19}$$

We are going to show that there exists a $C > 0$ such that

$$B \leq C\Bigg[\frac{1}{\sqrt{n}} + \sqrt{\frac{B}{n^{k/2}} + \sigma^2}\Bigg(\log\frac{A\|F\|_{L_2(P^m)}}{\sigma}\Bigg)^{k/2}\Bigg]. \tag{3.20}$$

We shall consider two cases. In case 1, $E\|U_n^{(k)}((P^{m-k}f)^2)\| \leq \sigma^2$. In this case, since the function $\sqrt{x}(-\log x)^{k/2}$, $0 < x < 1$, is increasing on $(0, e^{-k}]$, we get the trivial bound from (3.15) and (3.19) that

$$B \leq \frac{D}{\sqrt{n}} + 2\sigma\Bigg(\log\frac{2A\|F\|_{L_2(P^m)}}{\sigma}\Bigg)^{k/2},$$

which of course implies the bound (3.20) (note that $A\|F\|_{L_2(P^m)}/\sigma \geq e^m$).

Next consider case 2, $E\|U_n^{(k)}((P^{m-k}f)^2)\| > \sigma^2$. To handle this case we must bound $E\|U_n^{(k)}((P^{m-k}f)^2)\|$. It is here that we will use the induction hypothesis. By Hoeffding's decomposition (3.3) we have

$$E\|U_n^{(k)}((P^{m-k}f)^2)\| \leq \sum_{r=0}^{k}\binom{k}{r}E\|U_n^{(r)}(\pi_r(P^{m-k}f)^2)\|. \tag{3.21}$$

The term corresponding to $r = 0$ is simply $\|P^k(P^{m-k}f)^2\| \leq \sigma^2$. Consider now the new class

$$\mathcal{G} := \{(P^{m-k}f)^2 : f \in \mathcal{F}\}$$

of functions of $k$ variables. For any probability measure $Q$ on $S^k$,

$$Q((P^{m-k}f)^2 - (P^{m-k}g)^2)^2 \leq 4Q \times P^{m-k}(f-g)^2,$$

and thus the class $\mathcal{G}$ is VC-type with constants $A$ and $v$ and envelope $2\sqrt{P^{m-k}F^2}$ as in (2.9). Also observe that since

$$\|P^k(P^{m-k}f)^4\| \leq \|P^m f^2\| \leq 4\sigma^2 \leq 4\|\sqrt{P^{m-k}F^2}\|_{L_2(P^k)}^2 = 4\|F\|_{L_2(P^m)}^2,$$

we verify that (3.11) holds. Also we trivially see that (3.10) is satisfied, that is,

$$4n\sigma^2 > n\sigma^2 \geq C_2 \log\left(\frac{2\|F\|_{L_2(P^m)}}{\sigma}\right) = C_2 \log\left(\frac{4\|F\|_{L_2(P^m)}}{2\sigma}\right).$$

Finally noting that $A \geq e^m > e^k$ and $2\sqrt{P^{m-k}F^2} \leq 2$, we are permitted to apply the induction hypothesis for $1 \leq r < k$ ($r = 0$ has already been dealt with) to get with $M = 2$,

$$E\|U_n^{(r)}(\pi_r(P^{m-k}f)^2)\| \leq C_1(k, A, v, 2)n^{-r/2} 2\sigma\left(\log \frac{A\|F\|_{L_2(P^m)}}{\sigma}\right)^{r/2}. \tag{3.22}$$

Since by the hypotheses (3.10) and (3.11) on $\sigma$ [note (3.11) and $F \leq 1$ imply $\sigma^2 \leq 1$] we have

$$\frac{1}{n^{r/2}}\sigma\left(\log \frac{A\|F\|_{L_2(P^m)}}{\sigma}\right)^{r/2} \leq C^{r/2}\sigma^{1+r} \leq C^{r/2}\sigma^2, \tag{3.23}$$

it follows from the bounds (3.22) and (3.23) that

$$\binom{k}{r} E\|U_n^{(r)}(\pi_r(P^{m-k}f)^2)\| \leq C\sigma^2, \qquad 1 \leq r < k. \tag{3.24}$$

As for $r = k$, we randomize and use the entropy bound for Rademacher chaos, just as we did at the beginning of the proof, using the fact that $\mathcal{G}$ is VC-type. This gives

$$\begin{aligned} &E\|U_n^{(k)}(\pi_k(P^{m-k}f)^2)\| \\ &\quad \leq C|I_n^k|^{-1/2} E\left\| |I_n^k|^{-1/2} \sum_{I_n^k} \varepsilon_{i_1}\cdots\varepsilon_{i_k}(P^{m-k}f)^2 \right\| \\ &\quad \leq \frac{C}{|I_n^k|^{1/2}} E\int_0^{\sqrt{\|U_n^{(k)}((P^{m-k}f)^4)\|}} \left(\log \frac{A\sqrt{U_n^{(k)}(P^{m-k}F^2)}}{\tau}\right)^{k/2} d\tau \\ &\quad \leq C\frac{B}{n^{k/2}}, \end{aligned} \tag{3.25}$$

where $B$ is as defined in (3.13) and we use

$$\|U_n^{(k)}((P^{m-k}f)^4)\| \leq \|U_n^{(k)}((P^{m-k}f)^2)\|.$$

From (3.15), (3.19), (3.21)–(3.25) and $E\|U_n^{(k)}((P^{m-k}f)^2)\| > \sigma^2$ we get (3.20).

Since $\log[A\|F\|_{L_2(P^m)}/\sigma] > 1$, hypothesis (3.10) on $\sigma$ gives

$$\frac{1}{\sqrt{n}} \leq C\frac{\sigma}{(\log A\|F\|_{L_2(P^m)}/\sigma)^{1/2}} \leq C\sigma\left(\log \frac{A\|F\|_{L_2(P^m)}}{\sigma}\right)^{k/2},$$

and therefore it follows from inequality (3.20) that with perhaps a different value of $C$,

$$B \le C\sqrt{\frac{B}{n^{k/2}} + \sigma^2}\left(\log \frac{A\|F\|_{L_2(P^m)}}{\sigma}\right)^{k/2}.$$

Taking squares and solving this inequality, we get that there exists a constant $C$ (easy to evaluate) such that

$$B \le C\left[\frac{(\log A\|F\|_{L_2(P^m)}/\sigma)^k}{n^{k/2}} + \sigma\left(\log \frac{A\|F\|_{L_2(P^m)}}{\sigma}\right)^{k/2}\right].$$

But, again by condition (3.10) on $\sigma$ as in (3.23),

$$\frac{(\log A\|F\|_{L_2(P^m)}/\sigma)^{k/2}}{n^{k/2}} \le C\sigma^k \le \sigma,$$

and therefore

$$B \le C\sigma\left(\log \frac{A\|F\|_{L_2(P^m)}}{\sigma}\right)^{k/2},$$

which, by inequality (3.12), proves the theorem. □

REMARK 4. An analogue of the previous theorem holds if we replace the function $(\log(A\|F\|_{L_2(Q)}/\tau)$ by $H(\|F\|_{L_2(Q)}/\tau)$ with $H$ an increasing regularly varying function of exponent $0 \le \alpha < 2/m$, very much as in Theorem 3.1 in Giné and Koltchinskii [17].

3.4. *$U$-statistics taking values in separable Banach spaces.*

3.4.1. *The $L_p$ case,* $2 \le p < \infty$. We shall begin by establishing an exponential bound for the tail of the norm of a $U$-statistic taking values in a separable type 2 Banach space $\mathbf{B}$. Let us recall that a separable Banach space $\mathbf{B}$ is of type 2 if for any finite number of points $x_i \in \mathbf{B}$ and independent Rademacher variables $\varepsilon_i$ (independent random signs), we have $E\|\sum_{i=1}^n \varepsilon_i x_i\|^2 \le C\sum_{i=1}^n \|x_i\|^2$, and the smallest such constant $C$ is the type 2 constant of $\mathbf{B}$. It is well known that the $L_p$ spaces are of type 2 if (and only if) $p \ge 2$. If $Z_i$ are independent, centered $\mathbf{B}$-valued random vectors with a square integrable norm and $\mathbf{B}$ is of type 2, then

$$E\left\|\sum_{i=1}^n Z_i\right\|^2 \le C\sum_{i=1}^n E\|Z_i\|^2 \tag{3.26}$$

(Araujo and Giné [2] or Ledoux and Talagrand [24]).

THEOREM 9. *Let $H$ be a function of $k$ variables, $P$-canonical, symmetric, with values in a type 2 Banach space. We have, for all $x \geq 0$,*

$$\Pr\left\{\left\|\sum_{\mathbf{i}\in I_n^k} H(X_{i_1},\dots,X_{i_k})\right\| \geq x\right\} \leq D\exp\left\{-\left(\frac{x}{\lambda_0 n^{k/2}\kappa}\right)^{2/k}\right\}, \tag{3.27}$$

*where $\kappa$ is a number satisfying*

$$\kappa \geq \sup_{x_1,\dots,x_k\in S} \|H(x_1,\dots,x_k)\| \tag{3.28}$$

*and $\lambda_0$ and $D$ are constants that depend only on $k$ and the type 2 constant $C$ of $\mathbf{B}$.*

PROOF. (This inequality was mentioned on page 252 of de la Peña and Giné [7] and its precise statement and proof were left to the reader.) Here is the proof.

Let

$$S_n = \frac{1}{\kappa\binom{n}{k}^{1/2}} \sum_{1\leq i_1<\cdots<i_k\leq n} H(X_{i_1},\dots,X_{i_k}),$$

and let

$$S_n' = \frac{1}{\kappa\binom{n}{k}^{1/2}} \sum_{1\leq i_1<\cdots<i_k\leq n} \varepsilon_{i_1}^1\cdots\varepsilon_{i_k}^k H(X_{i_1}^{(1)},\dots,X_{i_k}^{(k)}),$$

where $\kappa$ is as in (3.28), $X_i^{(j)}$ are i.i.d. with law $P$, and $\varepsilon_i^j$ are i.i.d. Rademacher variables independent of the collection of variables $\{X_i^{(j)}\}$.

By decoupling (Theorem 3.1.1 of de la Peña and Giné [7]) and convexity (e.g., Theorem 3.5.3 on page 140 of de la Peña and Giné [7], as it applies to a nonnegative, nondecreasing convex function $\Psi$), there is a constant $C_k$ such that

$$\|\|S_n\|\|_{\Psi_{2/k}} \leq C_k \|\|S_n'\|\|_{\Psi_{2/k}}. \tag{3.29}$$

By hypercontractivity of Rademacher chaos (Khinchin's inequality for Rademacher chaos), for example, Theorem 3.2.2 on page 113 of de la Peña and Giné **(year?)**, for all $r \geq 2$,

$$E_\varepsilon\|S_n'\|^r \leq r^{kr/2}(E_\varepsilon\|S_n'\|^2)^{r/2},$$

and by the type 2 inequality applied one sequence $\{\varepsilon_i^j\}$, $j = 1,\dots,k$, at a time, we get

$$(E_\varepsilon\|S_n'\|^2)^{r/2} \leq C^{kr/2},$$

where $C$ is the type 2 constant of the Banach space, and therefore,

$$E_\varepsilon\|S_n'\|^r \leq C^{kr/2} r^{kr/2}.$$

This inequality yields, by Taylor expansion of the exponential, that for some $c > 0$,

$$E_\varepsilon \Psi_{2/k}(c\|S_n'\|) < \infty$$

and therefore that there exists a constant $\lambda_1 = \lambda_1(C,k)$ depending only on $C$ and $k$, such that

$$E_\varepsilon \Psi_{2/k}(\|S_n'\|/\lambda_1) \le 1.$$

Integrating with respect to the $X$'s,

$$E\Psi_{2/k}(\|S_n'\|/\lambda_1) \le 1, \qquad \text{that is, } \|\|S_n'\|\|_{\Psi_{2/k}} \le \lambda_1,$$

which, combined with (3.29), gives

$$\|\|S_n\|\|_{\Psi_{2/k}} \le C_k \lambda_1 := \lambda_0,$$

a constant. Of course this implies, by the definition of the Orlicz norm $\Psi_{2/k}$, that there is a constant $D(k,C)$ such that

$$\Pr\{\|S_n\| > x\} \le D e^{-(x/\lambda_0)^{2/k}}. \qquad \square$$

3.4.2. *The $L_p$ case*, $1 \le p < 2$. Our aim now is to obtain a useful exponential inequality for $\mathbf{B}$-valued $U$-statistics, when $\mathbf{B}$ is not necessarily a type 2 Banach space. We will generalize to $\mathbf{B}$-valued $U$-statistics (via decoupling followed by iteration) the following sharp inequality for sums of independent random vectors: there is a constant $L < \infty$ such that if $\mathbf{B}$ is a separable Banach space, $Z_i$, $i \in \mathbf{N}$, are independent mean zero random vectors taking values in $\mathbf{B}$ and $r \ge 2$, then, setting $S_n = \sum_{i=1}^n Z_i$,

$$E\|S_n\|^r \le L^r \left[ r^{r/2}(E\|S_n\|^2)^{r/2} + r^r E \max_{1\le i\le n} \|Z_i\|^r \right]. \tag{3.30}$$

This inequality was obtained by Pinelis [28], and it also follows easily from the sharper inequality in Giné, Latała and Zinn [18], Proposition 3.1. Next we extend inequality (3.30) to $\mathbf{B}$-valued $U$-statistics.

THEOREM 10. *Let $\mathbf{B}$ be a separable Banach space, let $H : S^k \mapsto \mathbf{B}$ be a bounded $P$-canonical random vector symmetric in its entries and let $X_i, X_i^{(j)}$ be i.i.d. $S$-valued random variables, $1 \le i \le n$ and $j = 1, \ldots, k$. Define $\kappa$ and $\chi_n$ to be any pair of numbers such that*

$$\begin{aligned} &\kappa \ge \sup_{x_1,\ldots,x_k \in S} \|H(x_1,\ldots,x_k)\| \\ &\textit{and} \quad \chi_n \ge \left( E \left\| \sum_{\mathbf{i}\in I_n^k} H(X_{i_1}^{(1)},\ldots,X_{i_k}^{(k)}) \right\|^2 \right)^{1/2}. \end{aligned} \tag{3.31}$$

*Then there exists a constant $C$ depending only on $k$ such that, for all $n \in \mathbf{N}$ and $r \geq 2$,*

$$E\left\|\sum_{\mathbf{i}\in I_n^k} H(X_{i_1}^{(1)},\ldots,X_{i_k}^{(k)})\right\|^r \leq C^r[r^{kr/2}\chi_n^r + r^{(k+1)r/2}n^{(k-1)r}\kappa^r + r^{kr}\kappa^r]. \tag{3.32}$$

*Moreover, there exists a constant $D_0$ depending only on $k$ such that, for all $n \in \mathbf{N}$ and $r \geq 2$,*

$$E\left\|\sum_{\mathbf{i}\in I_n^k} H(X_{i_1},\ldots,X_{i_k})\right\|^r \leq D_0^r[r^{kr/2}\chi_n^r + r^{(k+1)r/2}n^{(k-1)r}\kappa^r + r^{kr}\kappa^r], \tag{3.33}$$

*where $\kappa$ and $\chi_n$ are defined as in (3.31).*

PROOF. Inequality (3.30) gives the result for $k = 1$. Assume the result is true for $k$ and for every Banach space, and let $H$ be a function of $k+1$ variables satisfying the conditions in the statement of the theorem. Before starting the induction, we describe some simplifying notation. We will denote by $\leq'$ inequality up to a multiplicative constant $C^r$ with $C$ depending only on $k$. Also, we will write $H_{\mathbf{i}} = H$ if the coordinates of the multi-index $\mathbf{i} \in \{1,\ldots,k+1\}^n$ are all different, and $H_{\mathbf{i}} = 0$ otherwise, and we will drop the arguments of $H$ or $H_{\mathbf{i}}$, so $\sum_{\mathbf{i}} H_{\mathbf{i}}$ will mean $\sum_{\mathbf{i}\in\{1,\ldots,k+1\}^n} H_{\mathbf{i}}(X_{i_1}^{(1)},\ldots,X_{i_{k+1}}^{(k+1)})$. Finally, for $A \subset \{1,\ldots,k+1\}$, $E_A$ will mean integration with respect to the variables $X_i^{(r)}$ only, for $r \in A$ and $i \leq n$.

Applying (3.30) conditionally on the variables $X_i^{(j)}$, $j = 1,\ldots,k$, we obtain

$$\begin{aligned} E\left\|\sum_{\mathbf{i}} H_{\mathbf{i}}\right\| &\leq' E_{1,\ldots,k}\left[r^{r/2}\left(E_{k+1}\left\|\sum_{i_{k+1}}\left(\sum_{i_1,\ldots,i_k} H_{\mathbf{i}}\right)\right\|^2\right)^{r/2}\right] \\ &\quad + E_{1,\ldots,k}\left[r^r E_{k+1}\max_{i_{k+1}}\left\|\sum_{i_1,\ldots,i_k} H_{\mathbf{i}}\right\|^r\right]. \end{aligned} \tag{3.34}$$

In order to deal with the first term in (3.34), for each $x_1,\ldots,x_k \in S$ fixed we consider the random variable $\sum_{i_{k+1}} H_{\mathbf{i}}(x_1,\ldots,x_k,X_{i_{k+1}}^{(k+1)})$ as a function from $S^k$ into the Banach space $L_2(\Omega,\Sigma,P;\mathbf{B})$ of $\mathbf{B}$-valued random variables whose $\mathbf{B}$ norms are square integrable, with norm

$$|F(x_1,\ldots,x_k,X_{i_{k+1}}^{(k+1)})|_* := (E_{k+1}\|F(x_1,\ldots,x_k,X_{i_{k+1}}^{(k+1)})\|^2)^{1/2}.$$

To apply the induction hypothesis to the statistic $\sum_{i_1,\ldots,i_k}(\sum_{i_{k+1}} H_{\mathbf{i}})$ in this Banach space with the norm $|\cdot|_*$, we first note that, if we denote by $\tilde{\kappa}$ and

$\tilde{\chi}$ the corresponding quantities associated with this statistic, then

$$\tilde{\kappa} = \max_{x_1,\dots,x_k\in S}\left|\sum_{i_{k+1}} H(x_1,\dots,x_k,X^{(k+1)}_{i_{k+1}})\right|_* \le n\kappa$$

and

$$\tilde{\chi}_n = \left(E_{1,\dots,k}\left|\sum_{i_1,\dots,i_k}\left(\sum_{i_{k+1}} H_{\mathbf{i}}\right)\right|_*^2\right)^{1/2} = \left(E\left\|\sum_{\mathbf{i}} H_{\mathbf{i}}\right\|^2\right)^{1/2} = \chi_n.$$

Hence, the induction hypothesis gives

$$\begin{aligned} & r^{r/2}E_{1,\dots,k}\left(E_{k+1}\left\|\sum_{i_{k+1}}\left(\sum_{i_1,\dots,i_k} H_{\mathbf{i}}\right)\right\|^2\right)^{r/2} \\ &\quad = r^{r/2}E_{1,\dots,k}\left|\sum_{i_1,\dots,i_k}\left(\sum_{i_{k+1}} H_{\mathbf{i}}\right)\right|^r \\ &\quad \le' r^{r/2}[r^{kr/2}\tilde{\chi}_n^r + r^{(k+1)r/2}n^{(k-1)r}\tilde{\kappa}^r + r^{kr}\tilde{\kappa}^r] \\ &\quad \le r^{(k+1)r/2}\chi_n^r + r^{(k+2)r/2}n^{kr}\kappa^r + r^{(k+1/2)r}n^r\kappa^r. \end{aligned}$$

Now, if $n \le r^{1/2}$, then $r^{(k+1/2)r}n^r \le r^{(k+1)r}$ and if $n > r^{1/2}$, then $r^{(k+1/2)r}n^r \le r^{(k+2)r/2}n^{kr}$, so that

$$\begin{aligned} & r^{r/2}E_{1,\dots,k}\left(E_{k+1}\left\|\sum_{i_{k+1}}\left(\sum_{i_1,\dots,i_k} H_{\mathbf{i}}\right)^2\right\|\right)^{r/2} \\ &\quad \le' r^{(k+1)r/2}\chi_n^r + r^{(k+2)r/2}n^{kr}\kappa^r + r^{(k+1)r}\kappa^r. \end{aligned} \tag{3.35}$$

Finally, we apply the induction hypothesis conditionally on the variables $X_i^{(k+1)}$ to the second term of (3.34) by considering that term to be a $U$-statistic with values in the Banach space $\ell_n^\infty(\mathbf{B}) := \{(v_1,\dots,v_n) : v_i \in \mathbf{B}\}$, with norm $|(v_1,\dots,v_n)| := \max_{1\le i\le n}\|v_i\|$. In fact, with this definition,

$$\begin{aligned} & \max_{i_{k+1}}\left\|\sum_{i_1,\dots,i_k} H_{\mathbf{i}}\right\| \\ &\quad = \left|\sum_{i_1,\dots,i_k}(H_{\mathbf{i}}(X^{(1)}_{i_1},\dots,X^{(k)}_{i_k},X^{(k+1)}_1),\dots,H_{\mathbf{i}}(X^{(1)}_{i_1},\dots,X^{(k)}_{i_k},X^{(k+1)}_n))\right|, \end{aligned}$$

and if we denote by $\overline{\kappa}$ and $\overline{\chi}_n$ the corresponding parameters, we have, for each value of $X_i^{(k+1)}$, $j = 1,\dots,n$,

$$\overline{\kappa} = \max_{x_1,\dots,x_k\in S}\|(H(x_1,\dots,x_k,X^{(k+1)}_1),\dots,H(x_1,\dots,x_k,X^{(k+1)}_n))\| \le \kappa$$

and

$$\overline{\chi}_n = \left( E_{1,\dots,k} \max_{i_{k+1}} \left\| \sum_{i_1,\dots,i_k} H_{\mathbf{i}} \right\|^2 \right)^{1/2} \le n^k \kappa.$$

So, induction gives

$$r^r E_{1,\dots,k} \max_{i_{k+1}} \left\| \sum_{i_1,\dots,i_k} H_{\mathbf{i}} \right\|^r \le' r^r [r^{kr/2} \overline{\chi}_n^r + r^{(k+1)r/2} n^{(k-1)r} \overline{\kappa}^r + r^{kr} \overline{\kappa}^r]$$

$$\le \kappa^r [r^{(k+2)r/2} n^{kr} + r^{(k+3)r/2} n^{(k-1)r} + r^{(k+1)r}].$$

By considering the cases $n \le r^{1/2}$ and $n > r^{1/2}$ we see that

$$r^{(k+3)r/2} n^{(k-1)r} \le \max(r^{(k+2)r/2} n^{kr}, r^{(k+1)r}),$$

from which it follows that

$$r^r E_{1,\dots,k} \max_{i_{k+1}} \left\| \sum_{i_1,\dots,i_k} H_{\mathbf{i}} \right\|^r \le' \kappa^r [r^{(k+2)r/2} n^{kr} + r^{(k+1)r}]. \tag{3.36}$$

Now, the first part of the theorem follows by substituting the estimates (3.35) and (3.36) into inequality (3.34). The proof of the theorem is completed by noting that (3.33) follows from (3.32) via de la Peña's [6] decoupling inequality (see Theorem 3.1.1 in de la Peña and Giné [7]). □

We note that this theorem could have been proved, with only formal changes, for $H_{\mathbf{i}}$ depending on the subindices $\mathbf{i}$. Giné, Latała and Zinn, unpublished, have a moment inequality for $\mathbf{B}$-valued canonical $U$-statistics of order 2 without boundedness assumptions that contains (3.32) for $k = 2$. Our proof is inspired by theirs.

We shall apply (3.33) to the special case when $\mathbf{B} = L_p(\mathbf{R}^d)$, $1 \le p < 2$, and $H(x_1, \dots, x_k) = \pi_k \overline{K}_h(\cdot, x_1, \dots, x_k)$, to obtain the following. Recall the definition (2.21) of $\mu_s$.

COROLLARY 3. *Let $K$, $f_g$ and $s > 0$ satisfy the conditions of Theorem 3 for $1 \le p < 2$. Then there exist $\gamma > 0$ and $D > 0$ such for all $n \ge 1$ and $x > 0$,*

$$\Pr\left\{ \frac{\| \sum_{I_n^k} \pi_k \overline{K}_h(\cdot, X_{i_1}, \dots, X_{i_k}) \|_p}{D n^{k-1}/\sqrt{h}} > x \right\} \le \gamma \exp\left( -\left( \frac{x}{\|K\|_{p,2,s}} \right)^{1/k} \right), \tag{3.37}$$

*where $\|K\|_{p,2,s} = \|K\|_p \vee \|K\|_{L_2(\mu_s)}$.*

For the proof, first, note the following easy bound for $\|\|\pi_k \overline{K}_h\|_p\|_\infty$, namely, we can choose some $A_1 > 0$, so that

$$\kappa := A_1 \|K\|_p / h^{1-1/p} \ge \sup_{x_1,\dots,x_k \in S} \|\pi_k \overline{K}_h(\cdot, x_1, \dots, x_k)\|_p. \tag{3.38}$$

To see this, just observe that $\pi_k \overline{K}_h(t, x_1, \dots, x_k)$ is a linear combination of terms of the form $h_\ell(t, x_{i_1}, \dots, x_{i_\ell})$ with $(i_1, \dots, i_\ell) \in I_n^\ell$ and $0 \le \ell \le k$, where $h_0(t) = P^m K_h$ and, for $\ell \ge 1$,

$$\begin{aligned} h_\ell(t, x_{i_1}, \dots, x_{i_\ell}) &= P^{m-\ell} \overline{K}_h(t, x_1, \dots, x_\ell, X_{\ell+1}, \dots, X_m) \\ &= P^{m-\ell} \frac{1}{m!} \sum_{\sigma \in I_m^m} K_h(t - g_\sigma(x_1, \dots, x_\ell, X_{\ell+1}, \dots, X_m)) \end{aligned}$$

and $g_\sigma(y_1, \dots, y_m) = g(y_{\sigma_1}, \dots, y_{\sigma_m})$. Then Jensen's inequality for expectations and averages and a substitution yield

$$\int |h_\ell(t, x_1, \dots, x_l)|^p \, dt \le \int |K_h(t)|^p \, dt = h^{1-p} \int |K(t)|^p \, dt.$$

We shall also need a good bound for $E\| \sum_{I_n^k} \pi_k \overline{K}_h(\cdot, X_{i_1}, \dots, X_{i_k})\|_p$ in order to estimate $\chi_n$ in (3.31). Such a bound will be based on the following lemma, which is an extension of ideas and results in Chapter 7 in Devroye [8] and Section 3 of Devroye [9].

LEMMA 1. *Let $1 \le p < 2$. If $f$ and $k$ are two nonnegative functions on $\mathbf{R}^d$ such that $f, k \in L_1(\mu_s)$ for some $s > d(2-p)/p$, then, for any $0 < b < \infty$,*

$$\sup_{h \in (0,b]} \int (k_h * f)^{p/2}(y) \, dy \le C(\|f\|_{L_1(\mu_s)} \|k\|_{L_1(\mu_s)})^{p/2}$$

*for some constant $C$ that depends only on $d$, $p$ and $s$.*

PROOF. Let $q = sp/(2-p)$ and note that $v(y) := (1 + |y|)^{-q}$ is integrable. Then, by Jensen with respect to the probability measure $v(y)\,dy/\|v\|_1$, we have

$$\begin{aligned} \int (k_h * f)^{p/2}(y) \, dy &= \int v(y)^{-1} (k_h * f)^{p/2}(y) v(y) \, dy \\ &\le \|v\|_1^{1-p/2} \left( \int (1 + |y|)^s (k_h * f)(y) \, dy \right)^{p/2}. \end{aligned}$$

Since $(1 + |u + v|) \le (1 + |u|)(1 + |v|)$ and $1 + |hu| \le (1 \vee h)(1 + |u|)$ for $h > 0$, we also have

$$\begin{aligned} &\int (1 + |y|)^s (k_h * f)(y) \, dy \\ &\quad \le \int \int (1 + |y - x|)^s (1 + |x|)^s k_h(y - x) f(x) \, dx \, dy \\ &\quad \le (1 \vee b^{s/d}) \|f\|_{L_1(\mu_s)} \|k\|_{L_1(\mu_s)}. \end{aligned}$$

The lemma follows from these inequalities. □

LEMMA 2. *If $K \in L_2(\mathbf{R}^d)$ and $h > 0$, then, for $1 \le k \le m$ and $t \in \mathbf{R}^d$,*

$$E(\pi_k \overline{K}_h(t, X_1, \dots, X_k))^2 \le \frac{((K^2)_h * f_g)(t)}{h}. \tag{3.39}$$

*Let now $K$, $f_g$ and $s > 0$ satisfy the conditions of Theorem 3 for $1 \le p < 2$, and let $0 < b < \infty$. Then there exists $C < \infty$, such that, for $1 \le k \le m$ and all $0 < h \le b$,*

$$\left( E \left\| \sum_{I_n^k} \pi_k \overline{K}_h(\cdot, X_{i_1}, \dots, X_{i_k}) \right\|_p^2 \right)^{1/2} \le \frac{C \|K\|_{L_2(\mu_s)} n^{k/2}}{h^{1/2}}. \tag{3.40}$$

PROOF. Using the fact that $\pi_k$ is a projection in $L_2(\Pr)$, and applying convexity of the square function to the symmetrization of $K_h$, we have

$$\begin{aligned} E(\pi_k \overline{K}_h(t, X_1, \dots, X_k))^2 &\le \frac{1}{h^2} E K^2 \left( \frac{t - g(X_1, \dots, X_m)}{h^{1/d}} \right) \\ &= \frac{1}{h}((K^2)_h * f_g)(t), \end{aligned}$$

that is, (3.39). Next, by orthogonality,

$$E\left( \sum_{I_n^k} \pi_k \overline{K}_h(\cdot, X_{i_1}, \dots, X_{i_k}) \right)^2 = \frac{k! n!}{(n-k)!} E(\pi_k \overline{K}_h(t, X_1, \dots, X_k))^2.$$

Then the Minkowski inequality for integrals (e.g., Folland [14], page 194) and the above results yield

$$\begin{aligned} &\left( E \left\| \sum_{I_n^k} \pi_k \overline{K}_h(\cdot, X_{i_1}, \dots, X_{i_k}) \right\|_p^2 \right)^{1/2} \\ &\quad \le \left\| \left( E \left( \sum_{I_n^k} \pi_k \overline{K}_h(\cdot, X_{i_1}, \dots, X_{i_k}) \right)^2 \right)^{1/2} \right\|_p \\ &\quad \le \left( \frac{m! n^k}{h} \right)^{1/2} \left( \int ((K^2)_h * f_g)^{p/2}(t)\, dt \right)^{1/p}. \end{aligned}$$

The result (3.40) follows now from Lemma 1. □

PROOF OF COROLLARY 3. It follows from Lemma 2 and from (3.38) that we can take, for some $A_2 > 0$,

$$\chi_n := A_2 n^{k/2} \|K\|_{p,2,s} / \sqrt{h} \quad \text{and} \quad \kappa := A_1 \|K\|_{p,2,s} / h^{1-1/p}.$$

Assume that $0 < h \le 1$ and $r \ge 2$. Then we have the bounds

$$r^{k/2}\chi_n = A_2\|K\|_{p,2,s} r^{k/2} n^{k/2}/\sqrt{h} \le A_2\|K\|_{p,2,s} r^k n^{k-1}/\sqrt{h},$$

$$r^{(k+1)/2} n^{k-1}\kappa = A_1\|K\|_{p,2,s} r^{(k+1)/2} n^{k-1}/h^{1-1/p} \le A_1\|K\|_{p,2,s} r^k n^{k-1}/\sqrt{h}$$

and

$$r^k\kappa = r^k A_1\|K\|_{p,2,s}/h^{1-1/p} \le A_1\|K\|_{p,2,s} r^k n^{k-1}/\sqrt{h}.$$

By (3.33), this says that for some $D > 0$, for all $r \ge 2$,

$$E\left\|\sum_{I_n^k} \pi_k \overline{K}_h(\cdot, X_{i_1}, \dots, X_{i_k})\right\|_p^r \le (D\|K\|_{p,2,s} r^k n^{k-1}/\sqrt{h})^r$$

and thus

$$\begin{aligned} E\left\|\sum_{I_n^k} \pi_k \overline{K}_h(\cdot, X_{i_1}, \dots, X_{i_k})\right\|_p^{r/k} &\le \left(E\left\|\sum_{I_n^k} \pi_k \overline{K}_h(\cdot, X_{i_1}, \dots, X_{i_k})\right\|_p^r\right)^{1/k} \\ &\le (D\|K\|_{p,2,s} r^k n^{k-1}/\sqrt{h})^{r/k} \\ &= (3^k D\|K\|_{p,2,s} n^{k-1}/\sqrt{h})^{r/k} r^r/3^r. \end{aligned}$$

The same bound holds for $r = 1$. From this moment bound we easily get that for all $k \ge 2$,

$$E\exp\left(\left[\frac{\|\sum_{I_n^k} \pi_k \overline{K}_h(\cdot, X_{i_1}, \dots, X_{i_k})\|_p}{3^k D\|K\|_{p,2,s} n^{k-1}/\sqrt{h}}\right]^{1/k}\right) < \sum_{r=0}^{\infty} \frac{r^r}{3^r r!} =: \gamma < \infty.$$

The desired result (3.37) follows now from Markov's inequality and a renaming of $D$. □

Notice for future use that by (3.37) we get via (2.16) that for some $C > 0$,

$$\left\|\frac{\sqrt{n}\|\sum_{I_n^k} \pi_k \overline{K}_h(\cdot, X_{i_1}, \dots, X_{ik})\|_p}{n^k}\right\|_{\Psi_{1/k}} \le \frac{CD\|K\|_{p,2,s}}{\sqrt{nh}}. \tag{3.41}$$

3.5. *Completion of the proofs of the CLT.* We are now ready to prove (3.8) for $1 \le p \le \infty$, which will complete the proofs of Theorems 1, 2 and 3 (Section 3.1). We begin with the case $p = \infty$.

Proof of Theorem 1. Let

$$\mathcal{K}_{g,n} := \{K((y - g(\cdot))/h^{1/d}) : y \in \mathbf{R}^d, ah_n \le h \le bh_n\}.$$

Since $N(\mathcal{K}_{g,n}, L_2(\tilde{Q}), \varepsilon) = N(\mathcal{K}_n, L_2(\tilde{Q} \circ g^{-1}), \varepsilon)$ for every probability measure $\tilde{Q}$ on $\mathbf{R}^m$, it follows from condition (a) in Theorem 1 that the classes

$\mathcal{K}_{g,n}$ are VC-type with the same $A$ and $v$ and with envelopes $F_n(g)$. Moreover, since the map

$$(y, h, x_1, \dots, x_m) \mapsto \overline{K}((y - g(x_1, \dots, x_m))/h^{1/d})$$

is jointly measurable, these classes are image admissible Suslin, hence measurable. So, we can apply Theorem 8 to them. Since by Minkowski's inequality

$$\begin{aligned}(\lambda h_n)^2 E\overline{K}^2_{\lambda h_n}(y, X_1, \dots, X_m) &\le EK^2\left(\frac{y - g(X_1, \dots, X_m)}{(\lambda h_n)^{1/d}}\right) \\ &= \int K^2\left(\frac{y-w}{(\lambda h_n)^{1/d}}\right) f_g(w)\,dw \le \|f_g\|_\infty \|K\|_2^2 b h_n,\end{aligned}$$

we can take $\sigma^2 = Ch_n$ with $C = b\|f_g\|_\infty \|K\|_2^2 / \|\pi_k \overline{K}\|_\infty^2$. Hence, Theorem 8 gives that there is a constant $C_k$ such that, for all $n$ and $k$,

$$\begin{aligned}&E\left\{\sup\left|\frac{1}{n^{k/2}} \lambda h_n \sum_{I_n^k} \pi_k \overline{K}_{\lambda h_n}(t, X_{i_1}, \dots, X_{i_k})\right| : t \in \mathbf{R}^d, a \le \lambda \le b\right\} \\ &\qquad \le C_k \sqrt{h_n}[1 \vee \log(A\|F_n\|_{L_2(f_g)}/\sqrt{h_n}\,)].\end{aligned}$$

Hence, if $nh_n/[1 \vee \log(A\|F_n\|_{L_2(f_g)}/\sqrt{h_n}\,)]^2 \to \infty$, then for all $k = 2, \dots, m$,

$$E \sup_{a \le \lambda \le b} \|\sqrt{n} U_n^{(k)}(\pi_k \overline{K}_{\lambda h_n})\|_\infty \le \frac{C[1 \vee \log(A\|F_n\|_{L_2(f_g)}/\sqrt{h_n}\,)]}{n^{(k-1)/2} h_n^{1/2}} \to 0, \tag{3.42}$$

proving (3.8) in the $L_\infty$ case. □

Proof of Corollary 1. By the Hoeffding decomposition (3.5), it suffices to show [in analogy with Proposition 1 and (3.8)] that

$$\sup_{0<\lambda\le 1} |\lambda(\sqrt{n} m U_n^{(1)}(\pi_1 \overline{K}_{\lambda h_n}(0, \cdot, \dots, \cdot)) - \overline{\nu}_n(0))| \to 0 \qquad \text{in pr} \tag{3.43}$$

and

$$\sup_{0<\lambda\le 1} \lambda\sqrt{n}|U_n^{(k)}(\pi_k \overline{K}_{\lambda h_n}(0, \cdot, \dots, \cdot))| \to 0, \qquad k = 2, \dots, m. \tag{3.44}$$

As in the proof of Proposition 1, (3.43) reduces to proving

$$\sup_{0<\lambda\le 1} \lambda \int |\overline{\nu}_n(-u) - \overline{\nu}_n(0)||K_{\lambda h_n}(u)|\,du \to 0$$

for $i = 1, \dots, m$. Let $M$ be such that $D$ is contained in the ball of radius $M$ about the origin in $\mathbf{R}^d$. Since $K_{\lambda h_n}(u) = I\{u \in (\lambda h_n)^{1/d} D / (\lambda h_n \operatorname{Vol}(D))\}$ is zero for $|u| > (\lambda h_n)^{1/d} M$, we have

$$\sup_{0<\lambda\le 1} \lambda \int |\overline{\nu}_n(-u) - \overline{\nu}_n(0)||K_{\lambda h_n}(u)|\,du \le \sup_{|u| \le h_n^{1/d} M} |\overline{\nu}_n(-u) - \overline{\nu}_n(0)|,$$

which tends to zero in probability by condition (2.13), proving (3.43). For each $n$, the class of functions $\{I_{\lambda h_n D} : 0 \le \lambda \le 1\}$ is VC because it is linearly ordered (these are the simplest VC classes), and its envelope is contained in the class. The same is true for the class $\{\lambda h_n \overline{K}_{\lambda h_n}(0, \cdot, \dots, \cdot)/ \|\pi_k \overline{K}(0, \cdot, \dots, \cdot)\|_\infty : 0 < \lambda \le 1\}$. Then the logarithm in the bound (3.9) in Theorem 8 is simply a constant, and that theorem gives

$$n^{k/2} \frac{h_n}{\|\pi_k \overline{K}(0, \cdot, \dots, \cdot)\|_\infty} E \sup_{0<\lambda\le 1} |\lambda U_n^{(k)}(\pi_k \overline{K}_{\lambda h_n}(0, \cdot, \dots, \cdot))| \le C\sigma$$

as long as $n\sigma^2 > C'$, for fixed constants $C$ and $C'$ and for $k = 2, \dots, m$, where we can take

$$\sigma^2 = \frac{h_n \|f_g\|_\infty}{\mathrm{Vol}(D) \|\pi_k \overline{K}(0, \cdot, \dots, \cdot)\|_\infty^2}$$

because, by change of variables,

$$\sup_{0<\lambda\le 1} E[\lambda^2 h_n^2 K_{\lambda h_n}^2(-g(X_1, \dots, X_m)/(\lambda h_n)^{1/d})] \le h_n \|f_g\|_\infty / \mathrm{Vol}(D).$$

Since $nh_n \to \infty$ by hypothesis, the condition $n\sigma^2 > C'$ is satisfied, and therefore, we can apply Theorem 8 and obtain

$$\sup_{0<\lambda\le 1} \lambda\sqrt{n} |U_n^{(k)}(\pi_k \overline{K}_{\lambda h_n}(0, \cdot, \dots, \cdot))| \le \frac{C_k}{n^{(k-1)/2} h_n^{1/2}} \to 0$$

for all $2 \le k \le m$, proving (3.44). $\square$

To finish the proof of Theorems 2 and 3 is more complicated because we must deal with a mixed norm, the sup over $\lambda$ of an $L_p$ norm. The proof will consist in showing that the exponential inequalities from the above subsections lead to an entropy bound of the random variables in (3.8).

Given a kernel $L$ in $L_p(\mathbf{R}^d)$, $1 \le p < \infty$, and $h > 0$, let $\overline{L}_h(t, x_1, \dots, x_m)$ be the symmetrization of $L_h(t - g(x_1, \dots, x_m))$, as in (3.4). We observe that, just as in (3.38), there is a constant $c < \infty$ that depends only on $k$ such that

$$\| \|\pi_k \overline{L}_h\|_p \|_\infty \le \frac{c\|L\|_p}{h^{(p-1)/p}}, \qquad k = 1, \dots, m. \tag{3.45}$$

Because of (3.45), it makes sense to define

$$\mathbf{X}_{L,n} := \left\| \frac{h_n^{(p-1)/p}}{n^{k/2}} \sum_{\mathbf{i} \in I_n^k} \pi_k \overline{L}_{h_n}(\cdot, X_{i_1}, \dots, X_{i_m}) \right\|_p. \tag{3.46}$$

PROOF OF (3.8) FOR $2 \le p < \infty$. For $p \ge 2$, $L_p$ is of type 2 and Theorem 9 then gives that, for a constant $C$ that depends only on $k$ and $p$, for all $x \ge 0$,

$$\Pr\{\mathbf{X}_{L,n} \ge x\} \le D \exp\left\{ - \left( \frac{x}{C\|L\|_p} \right)^{2/k} \right\}, \tag{3.47}$$

and then, by (2.16), we have for another constant $C$ that depends on $k$ and $p$,

$$\|\mathbf{X}_{L,n}\|_{\Psi_{2/k}} \le C\|L\|_p, \tag{3.48}$$

where $\Psi_{2/k}$ is the Young modulus of exponential type defined in (2.14) [but with $\Psi_1(x) = e^x - 1$] and $\|\cdot\|_{\Psi_{2/k}}$ is the associated (pseudo)norm (2.15). Applying (3.48) to $L = K_\lambda - K_{\lambda'}$, we obtain

$$\|\mathbf{X}_{K_\lambda,n} - \mathbf{X}_{K_{\lambda'},n}\|_{\Psi_{2/k}} \le \|\mathbf{X}_{K_\lambda - K_{\lambda'},n}\|_{\Psi_{2/k}} \le C\|K_\lambda - K_{\lambda'}\|_p. \tag{3.49}$$

Then this bound and (3.45) allow us to apply the usual entropy integral bound, for example, in the version given in de la Peña and Giné [7], Corollary 5.1.5, and conclude that for some constant $D$, keeping in mind that $0 \in \mathcal{K}_{[a,b]}$,

$$\left\|\sup_{a\le\lambda\le b} \mathbf{X}_{K_\lambda,n}\right\|_{\Psi_{2/k}} \le D\int_0^{C\|K\|_p} \Psi_{2/k}^{-1}(N(\mathcal{K}_{[a,b]}, d_p, \varepsilon))\,d\varepsilon, \tag{3.50}$$

where $\mathcal{K}_{[a,b]}$ is defined by (2.17) and $d_p$ is the $L_p(\mathbf{R}^d)$ distance defined in (2.18) (technically, this only holds for any separable version of the process $\mathbf{X}_{K_\lambda,n}$, but we see in a remark below that this process itself is separable). Now, since up to constants $\Psi_\alpha$ is increasing in $\alpha$, hypothesis (a) in Theorem 2 implies that this integral is finite for all $2 \le k \le m$. Taking into account that the Orlicz distances of exponential type dominate (up to constants) the $L_r(\Pr)$ distances, inequality (3.50) implies that for all $2 \le k \le m$,

$$E \sup_{a\le\lambda\le b} \|\sqrt{n}U_n^{(k)}(\pi_k \overline{K}_{\lambda h_n})\|_p \le \frac{C}{n^{(k-1)/2} h_n^{(p-1)/p}},$$

for some constant $C$ that depends on $k, p, a$ and $b$. The condition $nh_n^{2(p-1)/p} \to \infty$ implies that this expectation tends to zero for $2 \le k \le m$, proving (3.8) for $p \ge 2$. $\square$

PROOF OF (3.8) FOR $1 \le p < 2$. Since $f_g$ and $K^2$ are in $L_1(\mu_s)$ for some $s > d(2-p)/p$, Lemma 1 applies and therefore Corollary 3 and inequality (3.41) hold. Now, the proof of (3.8) and subsequently Theorem 3 follow exactly as the proof of Theorem 2, but using Corollary 3 and its consequence (3.41) instead of Theorem 9 and its consequence for $\mathbf{X}_{L,n}$. These give

$$\|\mathbf{X}_{L,n}\|_{\Psi_{1/k}} \le Cn^{k/2-1}h_n^{1/2-1/p}\|L\|_{p,2,s}$$

instead of (3.48), which yields, as in the previous proof, for some constant $D$,

$$\left\|\sup_{a\le\lambda\le b} \mathbf{X}_{K_\lambda,n}\right\|_{\Psi_{1/k}}$$
$$\le Dn^{k/2-1}h_n^{1/2-1/p}\int_0^{C\|K\|_{p,2,s}} \Psi_{1/k}^{-1}(N(\mathcal{K}_{[a,b]}, d_p \vee \tilde{d}_{2,s}, \varepsilon))\,d\varepsilon$$

instead of (3.50). Hence, for some constant $C$,

$$E \sup_{a\le\lambda\le b} \|\sqrt{n}U_n^{(k)}(\pi_k\overline{K}_{\lambda h_n})\|_{L_p} \le' \frac{1}{n^{k/2-1/2}h_n^{1-1/p}} E \sup_{a\le\lambda\le b} \mathbf{X}_{K_\lambda,n}$$
$$\le \frac{C\|K\|_{p,2,s}}{(nh_n)^{1/2}} \to 0$$

if $nh_n \to \infty$. □

REMARK 5 (Separability of the process $\mathbf{X}_{K_\lambda,n}$). In the previous subsection technically we must make sure that

$$\left\|\sup_{\lambda\in[a,b]} \mathbf{X}_{K_\lambda,n}\right\|_{\Psi_\alpha} = \sup\left\{\left\|\max_{\lambda\in C} \mathbf{X}_{K_\lambda,n}\right\|_{\Psi_\alpha} : C \text{ finite, } C\subset[a,b]\right\}$$

in order to ensure that the entropy bound applies exactly to the process $\mathbf{X}_{K_\lambda,n}$ and not to a modification thereof. For this, by standard arguments (basically the monotone convergence theorem), it suffices that there exist $D\subset[a,b]$ countable such that

$$\sup_{\lambda\in[a,b]} \mathbf{X}_{K_\lambda,n} = \sup_{\lambda\in D} \mathbf{X}_{K_\lambda,n} \qquad \text{a.s.} \tag{3.51}$$

Let $D$ denote the rationals in $[a,b]$. To show (3.51) it suffices to prove that, with probability 1, for each $\lambda\in[a,b]$ and any sequence $\lambda_m$ in $D$ such that $\lambda_m\to\lambda$ we have $\lim_{m\to\infty}\mathbf{X}_{K_{\lambda_m}-K_\lambda,n}=0$. In turn to verify this it is enough to check that

$$\lim_{m\to\infty}\left(\int_{\mathbf{R}^d}\left|\sum_{I_n^k}\pi_k\overline{K_{\lambda_m h_n}-K_{\lambda h_n}}(t,X_{i_1},\dots,X_{i_k})\right|^p dt\right)^{1/p}=0. \tag{3.52}$$

Now, each of these $L_p$ norms is bounded by a finite number of terms of the form

$$\left(\int_{\mathbf{R}^d}|E(K_{\lambda_m h_n}-K_{\lambda h_n})(t-V)|^p\,dt\right)^{1/p},$$

where $V$ is a random variable. Observe that

$$\left(\int_{\mathbf{R}^d}|E(K_{\lambda_m h_n}-K_{\lambda h_n})(t-V)|^p\,dt\right)^{1/p}$$
$$\le\left(E\int_{\mathbf{R}^d}|(K_{\lambda_m h_n}-K_{\lambda h_n})(t-V)|^p\,dt\right)^{1/p}$$
$$=\left(E\int_{\mathbf{R}^d}\left|\frac{1}{\lambda_m h_n}K\left(\frac{t-V}{(\lambda_m h_n)^{1/d}}\right)-\frac{1}{\lambda h_n}K\left(\frac{t-V}{(\lambda h_n)^{1/d}}\right)\right|^p dt\right)^{1/p},$$

which by the change of variables inside the integral $u = (t - V)/(\lambda_m h_n)^{1/d}$ is equal to

$$\left(\lambda_m h_n \int_{\mathbf{R}^d} \left| \frac{1}{\lambda_m h_n} K(u) - \frac{1}{\lambda h_n} K\left(\left(\frac{\lambda_m}{\lambda}\right)^{1/d} u\right)\right|^p du\right)^{1/p}$$
$$= \left(\int_{\mathbf{R}^d} \left| K(u) - \frac{\lambda_m}{\lambda} K\left(\left(\frac{\lambda_m}{\lambda}\right)^{1/d} u\right)\right|^p du\right)^{1/p} (\lambda_m h_n)^{1/p-1}.$$

Now since $K \in L_p(\mathbf{R}^d)$ we have $\lim_{\gamma\to 1} \int_{\mathbf{R}^d} |K(u) - K(\gamma u)|^p \, du = 0$, and this, in turn, implies that

$$\lim_{m\to\infty} \int_{\mathbf{R}^d} \left| K(u) - \frac{\lambda_m}{\lambda} K\left(\left(\frac{\lambda_m}{\lambda}\right)^{1/d} u\right)\right|^p du = 0,$$

which gives (3.51).

**Acknowledgments.** We would like to thank Anton Schick and Wolfgang Wefelmeyer for introducing us to this area via the second named author. We are also grateful to an anonymous referee whose excellent detailed comments were very helpful for improving and streamlining the presentation; in particular the present proof of Lemma 1 belongs to him. Finally we thank Julia Dony for pointing out a number of misprints as well as for a helpful suggestion.

Department of Mathematics
University of Connecticut
196 Auditorium Road
Storrs, Connecticut 06269-3009
USA
E-mail: gine@math.uconn.edu

Department of Food and Resource Economics
University of Delaware
206 Townsend Hall
Newark, Delaware 19717
USA
E-mail: davidm@udel.edu